\documentclass[a4paper,12pt]{article} 
\newdimen\paperhight
\usepackage{amsmath}
\usepackage{amssymb}
\usepackage{amsfonts}
\usepackage[usenames]{color}

\newcommand{\dsp}{\displaystyle}

\newcommand{\hf}{\frac{1}{2}}

\newcommand{\pr}{\par \vspace{3mm}\noindent [{\bf Proof}] \qquad}
\newcommand{\prend}{\hfill \qed \par \vspace{3mm}}

\newcommand{\qed}{\quad\hbox{\rule[-2pt]{3pt}{6pt}}\par\vspace{3mm}}

\newcommand{\1}{{\bf 1}} 
 
\newcommand{\C}{\mathbb C} 
\newcommand{\Z}{\mathbb Z}

\newcommand{\N}{\mathbb N} 
 
\newcommand{\CC}{{\cal C}}

\newcommand{\CW}{{\cal W}}

\newcommand{\CH}{{\cal H}}

\newcommand{\FA}{{\mathfrak A}}
\newcommand{\FB}{{\mathfrak B}}
\newcommand{\FC}{{\mathfrak C}}
\newcommand{\FD}{{\mathfrak D}}
\newcommand{\FE}{{\mathfrak E}}
\newcommand{\FN}{{\mathfrak N}}
\newcommand{\al}{\alpha}
\newcommand{\be}{\beta}

\newcommand{\ga}{\gamma}
\newcommand{\la}{\lambda}
 
\newcommand{\tr}{{\rm tr}}

\newcommand{\wt}{{\rm wt}}
\newcommand{\Hom}{{\rm Hom}}
\newcommand{\Rad}{{\rm Rad}}
\newcommand{\soc}{{\rm soc}}
\newcommand{\Verma}{{\rm Verma}}
\newcommand{\End}{{\rm End}}
\newcommand{\Image}{{\rm Im}}
\newcommand{\Ker}{{\rm Ker}}
\newcommand{\Res}{{\rm Res}}
\newtheorem{thm}{Theorem}[section]
\newtheorem{prn}[thm]{Proposition}
\newtheorem{dfn}[thm]{Definition}
\newtheorem{cry}[thm]{Corollary}
\newtheorem{lmm}[thm]{Lemma}

\begin{document}

\title{Modular invariance of vertex operator 
algebras satisfying $C_{2}$-cofiniteness}
\author{Masahiko Miyamoto
\footnote{e-mail: miyamoto@math.tsukuba.ac.jp  \quad
Supported by the Grants-in-Aids for Scientific Research, 
No. 1344002, The Ministry of Education, Science and 
Culture, Japan.}} 
\date{\begin{tabular}{c}
Institute of Mathematics \cr
University of Tsukuba \cr
Tsukuba 305, Japan 
\end{tabular} }
\maketitle

\begin{abstract}
We investigate trace functions of modules for vertex operator 
algebras satisfying $C_2$-cofiniteness. For the modular invariance 
property, Zhu assumed two conditions in \cite{Zh}: $A(V)$ 
is semisimple and $C_2$-cofiniteness.  
We show that $C_2$-cofiniteness is enough to prove a modular 
invariance property. 
For example, if a VOA $V=\oplus_{m=0}^{\infty}V_m$ is $C_2$-cofinite, 
then the space spanned by generalized characters (pseudo-trace 
functions of the vacuum element) of $V$-modules is 
a finite dimensional $SL_2(\Z)$-invariant space and 
the central charge and conformal weights are all rational numbers. 
Namely we show that 
$C_2$-cofiniteness implies ``rational conformal field theory'' 
in a sense as expected in \cite{GN}. 
Viewing a trace map as one of symmetric linear maps and using a 
result of symmetric algebras,  
we introduce ``pseudo-traces'' and pseudo-trace functions and then 
show that the space spanned by 
such pseudo-trace functions has a modular invariance property. 
We also show that $C_2$-cofiniteness is equivalent to 
the condition that every weak module is an $\N$-graded weak 
module which is a direct sum of generalized 
eigenspaces of $L(0)$. 
\end{abstract}

\section{Introduction}
In this paper, we will consider a vertex operator algebra 
$V=\coprod_{n=0}^{\infty}V_n$ with central charge $c$. 
One of the central concepts in conformal 
field theory (CFT) is ``rationality,'' a condition 
which is supposed to 
express a kind of finiteness of the theory. There exists various 
notions of finiteness. One of them is the complete reducibility 
of $\N$-graded weak modules, which is the condition called 
rationality by the most researchers on vertex operator 
algebras. 
There is another important finiteness 
condition called ``$C_2$-cofiniteness''. 
Complete reduciblity of $\N$-graded weak modules is a 
condition for modules. In this case, 
it was shown in \cite{DLiM1} that Zhu algebra 
$A(V)=V/O(V)$ is a finite dimensional semisimple algebra and $V$ 
has only finitely many 
irreducible modules. On the other hand, $C_2$-cofiniteness 
is a property of $V$ itself, that is, 
if $C_2(V)=\left< a(\!-\!2)b \mid a,b\in V\right> $ 
is of finite codimension in $V$, then $V$ is called 
$C_2$-cofinite. Here 
$Y(a,z)=\sum_{n\in \Z}a(n)z^{-n-1}$ is the vertex 
operator of $a$. Then it was shown in \cite{DLiM3} that 
$A(V)$ is of finite dimension and so 
$V$ has only finitely many irreducible modules, 
say $\{W^1,...,W^k\}$.  
Under two hypotheses that $A(V)$ is semisimple and 
$V$ is $C_2$-cofinite, (Condition $C$ in his paper), 
Zhu proved in \cite{Zh} that the space spanned by trace 
functions 
$$ S^{W^i}(v, \tau)=\tr_{|W^i}o(v)q^{L(0)-c/24}
\qquad (q = e^{2\pi i \tau} \mbox{ and } \tau\in \CH)$$
has a modular-invariance 
($SL(2, \Z)$-invariance) property, where $c$ is the central 
charge of $V$, $\CH$ is the upper half plane 
$\{\tau\in \C \mid {\rm Im}(\tau)>0\}$ 
and $o(v)$ is the grade-preserving operator of $v\in V$. 
Namely, for $\binom{a\ b}{c\ d}\in SL(2,\Z)$, there is a 
$k\times k$-matrix $(\la_{ij})$ with $\la_{ij}\in \C$ 
such that 
$$\frac{1}{(c\tau\!+\!d)^n}
\left(S^{W^1}\!(u, \frac{a\tau\!+\!b}{c\tau\!+\!d})\cdots 
S^{W^k}\!(u, \frac{a\tau\!+\!b}{c\tau\!+\!d})\right)
= \left(S^{W^1}\!(u, \tau)\!\cdots \!S^{W^k}\!(u, \tau)\right)
\left(\begin{array}{ccc}\!\la_{11}&\!\!\cdots \!\!&\la_{1k}\!\cr 
\vdots&\!\!\cdots\!\!&\vdots \cr 
\!\la_{k1}&\!\!\cdots\!\!& \la_{kk}\!\end{array}\right) $$
for all $u\in V_{[n]}$ and $n=0,1,\cdots$, 
where $V=\oplus_{n=0}^{\infty}V_{[n]}$ is the 
second grading on $V$ introduced 
in \cite{Zh}, see Definition 2.3. 
Actually, Zhu assumed one more condition, but 
it is not necessary as mentioned in \cite{DLiM3}. 
This is a fundamental result for modular invariance 
properties of vertex operator algebras and 
is extended by several authors, e.g. 
\cite{DLiM3}, \cite{Miy1}-\cite{Miy3} and \cite{Y}. 

Since then, it becomes an important problem to study 
a relation between rationality and $C_2$-cofiniteness. 
It was once conjectured the equivalency between $C_2$-cofinitenss 
and rationality, but 
a few research don't suggest that $C_2$-cofiniteness implies 
the complete reducibility, see \cite{GK} and \cite{Mil}. 

On the other hand, $C_2$-cofiniteness was studied by Gaberdiel and 
Neitzke from a view point of rationality in \cite{GN} 
and they showed that if $V$ is a $C_2$-cofinite VOA of CFT type, 
then $V$ satisfies the most conditions required for a rational 
conformal field theory except for a modular invariance 
property, where $V$ is called CFT type 
if $V=\oplus_{n=0}^{\infty}V_n$ 
and $\dim V_0=1$, see also \cite{Zh}, \cite{DLiM3}, \cite{Li2}. 

Although $C_2$-cofiniteness was introduced by Zhu as a technical 
assumption and it is a property of VOA itself, it is a natural 
condition to consider the characters of all (weak) modules. 
For, in order to define $q^{L(0)}$ on a weak module $W$, 
$W$ has to be a direct sum of generalized eigenspaces of 
$L(0)$, which is a condition equivalent to $C_2$-cofiniteness. \\

\noindent
{\bf Theorem 2.7} \qquad 
{\it Let $V$ be a vertex operator algebra. 
Then the following are equivalent. \\
{\rm (1)} \quad $V$ is $C_2$-cofinite. \\
{\rm (2)} \quad Every weak module is a direct sum of 
generalized eigenspaces of $L(0)$. \\
{\rm (3)} \quad Every weak module is an 
$\N$-graded weak module $W=\oplus_{n=0}^{\infty}W(n)$ such that 
$W(n)$ is a direct 
sum of generalized eigenspaces of $L(0)$. \\
{\rm (4)} \quad $V$ is finitely generated and 
every weak module is an $\N$-graded weak module. }\\

\vspace{4mm}

Therefore if $V$ is $C_2$-cofinite, then 
we can define several modules naturally. For example, 
for $\N$-graded weak $V$-modules $U$, the maximal 
weak $V$-submodule $D(U)$ 
of $\Hom(U,\C)$ introduced in \cite{Li2} is an $\N$-graded weak 
module. The most results about $C_2$-cofininteness 
came from the existence of some spanning set, e.g. 
\cite{GN},\cite{Li2},\cite{Bu},\cite{ABD}, 
We will prove the existence of the following 
spanning set for a general VOA without negative weights. \\ 

\noindent
{\bf Lemma 2.4} \quad 
{\it Let $A$ be a set of homogeneous elements of $V$ 
such that $V=C_2(V)+\langle A\rangle$.  
Assume that $V$ is $C_2$-cofinite and 
$W$ is a weak module generated from 
$w$ (by the action of vertex operator). 
Then $W$ is spanned by 
the following elements  
$$  \{v^1(i_1)\cdots v^k(i_k)w   \ \mid \ v^i\in A, 
\quad i_1<\!\cdots \!<i_k \}.$$
In particular, if we set 
$$\CW(m)=\left< v^1(i_1)\cdots v^k(i_k)w \ \mid \ v^i\in V, \quad  
\deg(v^1(i_1)\cdots v^k(i_k))=m \right>,  $$
then $\CW(m)=0 \mbox{  for  } m\ll 0$, where 
$\deg(v(i)\cdots u(j))$ denotes the degree of $v(i)\cdots u(j)$ 
as an operator. }

\vspace{5mm}

The main purpose of this paper is to show that 
$C_{2}$-cofiniteness is enough for a modular invariance 
property without assuming that $A(V)$ is semisimple. 
As a result, the central charge $c$ and the conformal 
weights are all rational numbers.  Namely, $C_2$-cofiniteness 
provides the conditions required for a rational conformal 
field in a sense. 

In the proof of modular invariance property in \cite{Zh}, 
Zhu introduced the space $\CC_1(V)$ of 
one point functions, see also \cite{DLiM3}.  
He proved the modular invariance property by showing that
$\CC_1(V)$ is spanned by trace functions 
$S^W(v,\tau)=\tr_W o(v)q^{L(0)-c/24}$. 
However, if $A(V)$ is not semisimple, $\CC_1(V)$ may not be 
spanned by trace functions. 
One of our aims in this paper is to stuff suitable functions 
into a crevice. We introduce a new kind of trace map called 
a pseudo-trace map (different from pseudo trace in algebraic 
number field) on some kind of $\N$-graded weak $V$-modules with 
homogeneous spaces of finite dimension. 
One of the key steps in Zhu's proof is 
that $\CC_1(V)$ is spanned by functions with the form 
$$\sum_{i=0}^N(\sum_{j=1}^d(\sum_{k=0}^{\infty}
C_{ij,k}(v)q^k)q^{r_{ij}})(\ln(q))^i$$
such that the coefficient $C_{ij,0}$ of the lowest degree,
satisfies the conditions: \\
(1) $C_{ij,0}(O(V))=0$, \\
(2) $C_{ij,0}(ab\!-\!ba)=0$ for all $a,b\in A(V)$ and \\
(3) $C_{ij,0}((\omega\!-\!\frac{c}{24}\!-\!r_{ij})^{N-i+1}v)=0$. \\
In particular, $A(V)/{\rm Rad}(C_{ij,0})$ is a symmetric 
algebra with a symmetric linear function $C_{ij,0}$, where 
${\rm Rad}(\phi)\overset{{\rm def}}{=}
\{a\in A(V) \mid \phi(A(V)aA(V))=0\}$. 
We are not interested in symmetric algebras, but symmetric 
linear functions. 
Originally, Nesbitt and Scott showed in \cite{NS} 
that $A$ is a symmetric 
algebra if and only if its basic algebra $P$ is symmetric. 
This is an idea to explain our strategy in this paper. 
We will show that $C_{ij,n}$ is a symmetric linear map of 
$n$-th Zhu algebra $A_n(V)$ (see \S 2.2) and so 
we have a symmetric algebra 
$A=A_n(V)/\Rad(C_{ij,n})$ and 
its symmetric basic algebra $P$ with a symmetric linear 
function $\phi$. We start from $(P,\phi)$ and 
construct a right $P$-module $W$ (a generalized Verma module) 
such that the basic algebra of $R=\End_P(W)$ is $P$.  
We will call such a module $W$ ``interlocked with $\phi$.''
Nesbitt and Scott's result tells that $R$ has a symmetric 
linear map $\tr^{\phi}$ (we will call it {\bf pseudo-trace}). 
Then we will define a pseudo-trace function 
$$\tr^{\phi}_{W}o(v)q^{L(0)-c/24}$$ and show that $\CC_1(V)$ is 
spanned by such pseudo-trace functions. 

For example, 
$$P=\left\{ p=\left(\begin{array}{cc}
a &b\cr 0&a \end{array}\right)  \mid  
a,b\in \C \right\} $$
is a basic symmetric algebra with a linear map $\phi(p)=b$. 
We note $J(P)=\soc(P)=\left\{\left(\begin{array}{cc}
0&b\cr 0&0 \end{array}\right)  \mid b\in \C \right\}$. 
Consider a right $P$-module $T=\C^m\oplus \C^m$. 
Then 
$$ \End_P(T)=\left\{ \al=\left(\begin{array}{cc}
A_{\al}&B_{\al}\cr 0&A_{\al} \end{array}\right)  \mid  
A_{\al},B_{\al}\in M_{m,m}(\C) \right\} $$
and the basic algebra of $\End_P(T)$ is $P$. 
For any $\al\in \End_P(T)$, if we define 
$\tr^{\phi}_T(\al)=\tr \al:T/TJ(P)\to TJ$, that is,  
$\tr^{\phi}_T(\al)=\tr(B_{\al})$, then $\tr^{\phi}_T$ is 
also a symmetric linear map.  

We should note one more thing. Since we will 
treat the general cases, $L(0)$ may not act on $W^{(n)}_T(m)$ 
as a scalar. However, since $V$ is $C_2$-cofinite, every $n$-th 
Zhu algebra $A_n(V)$ is finite dimensional as we will see and 
so every generalized Verma $V$-module is 
a direct sum of modules $W=\oplus_{m=0}^{\infty}W(m)$ so that 
$L(0)\!-\!r\!-\!m$ acts on $W(m)$ as a nilpotent operator for 
some $r\in \C$, say $(L(0)\!-\!r\!-\!m)^s=0$ on $W(m)$. 
Let $L^s(0)$ denote the semisimple part of $L(0)$, that is, 
$m+r$ on $W(m)$. Then trace function $\tr^{\phi}_Wo(v)q^{L(0)}$ 
on $W$ is defined by 
$$\tr^{\phi}_Wo(v)q^{L(0)-c/24}=\tr^{\phi}_W\{o(v)
\sum_{i=0}^{s-1}\frac{(2\pi i\tau)^i(L(0)\!-\!L^s(0))^i}{i!}\}q^{L^s(0)-c/24}.$$ 
These modules $W$ are called ``logarithmic modules'', 
(see \cite{G},\cite{Mil}) and Flohr introduced in \cite{F} a 
concept of generalized characters to interpret the modular 
invariance property of characters of logarithmic modules. 
As we will see, $S^{W^{(n)}_T}(\1,\tau)$ is a linear 
combination of (ordinary) 
characters with coefficients in $\C[\tau]$ and plays a role of 
a generalized character. \\

Our main theorem is: \\

\noindent
{\bf Theorem 5.5}\qquad 
{\it Let $V$ be a $C_2$-cofinite VOA with central charge $c$. 
Then for $v\in V_{[m]}$, the set of pseudo-trace functions 
$$\left<    \tr^{\phi}_{W} o(v)q^{L(0)-c/24}   \mid  \mbox{ $W$ is 
interlocked with a symmetric linear 
map $\phi$ of $A_n(V)$} \right>  $$
is invariant under the action of $SL_2(\Z)$ with weight $m$. 
In particular, $\dim \CC_1(V)\!=\!
\dim A_n(V)/[A_n(V),A_n(V)]\!-\!
\dim A_{n-1}(V)/[A_{n-1}(V),A_{n-1}(V)]$ for $n\gg 0$. 
}\\

In particular, the space spanned by generalized 
characters $\tr^{\phi}_W q^{L(0)-c/24}$ is 
a finite dimensional $SL_2(\Z)$-invariant space. 
As a corollary, we obtain: \\

\noindent
{\bf Corollary 5.10 and 5.11}\qquad
{\it If $V$ is a $C_2$-cofinite VOA, 
then the central charge 
and the conformal weights are all rational numbers. 
Moreover, we have 
$$  \tilde{c}\leq \frac{\dim V/C_2(V)\!-\!1}{2},$$
where $\tilde{c}\overset{{\rm def}}{=}c\!-\!h_{min}$ is the 
effective central charge of $V$ 
and $h_{min}$ is the smallest conformal weight.} 
\\

They showed in \cite{GN} that 
$C_2$-cofiniteness implies $C_{m}$-cofiniteness for any 
$m=1,2,\cdots$ if $V$ is of CFT type. 
In this paper, we will consider only a VOA without negative weights. 
However, if we consider $C_{2+s}$-cofiniteness when 
$V=\oplus_{n=-s}^{\infty}V_n$ has a negative weight, 
then it is not difficult to see that we will have a similar result 
as in Lemma 2.4 and the other results in this paper by replacing 
$C_2(V)$ by $C_{2+s}(V)$ and $\wt(v)$ by $\wt(v)-s$, respectively. \\

This paper is organized as follows. 
In Section 2, we will explain the notation and fundamental 
results. In Section 3, we will introduce a concept of modules 
interlocked with $\phi$ and define a pseudo-trace map 
explicitly. In Section 4, we will define a pseudo-trace function 
for an $\N$-graded weak $V$-module interlocked with a symmetric 
function. In particular, we will explain that if we take a 
sufficiently large integer $n$ and $\phi$ is a symmetric 
linear map of $A_n(V)$, then we can construct a generalized 
Verma VOA-modules $W^{(n)}_T$, 
which is interlocked with $\phi$. In Section 5, we will prove a 
modular invariance property. \\

\noindent
{\bf Acknowledgement} \\
The author wishes to thank Hisaaki Fujita and Kenji Nishida for 
information about symmetric rings and Geoffrey Mason 
for information about irrational VOAs. 
He would also like to thank Matthias R. Gaberdiel for 
information about a generalized character.

\section{Fundamental results}
\subsection{Vertex operator algebras}

\begin{dfn}
A {\bf weak module} for VOA $(V, Y, {\bf 1}, \omega)$ is a 
vector space $M$; 
equipped with a formal power series 
$$Y^M(v, z)=\sum_{n \in {\bf Z}} v^M(n)z^{-n-1} 
\in ({\rm End}(M))[[z, z^{-1}]] $$ 
(called {\it the module vertex operator} of $v$) for $v \in V$ 
satisfying: \\
{\rm (1)}\quad  $v^M(n)w=0$ for $n\gg 0$ 
where $v\in V$ and $w\in M$; \\
{\rm (2)}\quad  $Y^M({\bf 1}, z)=1_M$; \\
{\rm (3)}\quad  $Y^M(\omega, z)
=\sum_{n\in\Z} L^M(n)z^{-n-2}$ satisfies: \\
$ \mbox{\ }$ \quad{\rm (3.a)}  the Virasoro algebra relations: 
$$ [L^M(n), L^M(m)]
=(n\!-\!m)L^M(n\!+\!m)+ \delta_{n+m, 0}\frac{n^3-n}{12}c, $$
$ \mbox{\ }$ \quad{\rm (3.b)}  the $L(-1)$-derivative property:
$$  Y^M(L(-1)v, z)=\frac{d}{dz}Y^M(v, z), $$
{\rm (4)}\quad  ``Commutativity'' holds; 
$$  [v^M(n), u^M(m)]=\sum_{i=0}^{\infty}\binom{n}{i}
\left(v(i)u\right)^M(n\!+\!m\!-\!i) \qquad \mbox{  and} $$
{\rm (5)}\quad  ``Associativity'' holds;  
$$ \left(v(n)u\right)^M(m)=\sum_{i=0}^{\infty}(-1)^i\binom{n}{i}
\left\{v^M(n\!-\!i)u^M(m\!+\!i)
-(\!-\!1)^nu^M(n\!+\!m\!-\!i)v^M(i) \right\},$$
where $\binom{n}{i}=\frac{n(n-1)\cdots (n-i+1)}{i!}$.\\

An {\bf $\N$-graded weak module} is a weak $V$-module which 
carries an ${\N}$-grading, 
$M=\oplus_{n=0}^{\infty}M(n)$, such that \\

\noindent
{\rm (1')}\quad  if $v\in V_r$ then $v^M(m)M(n)\subseteq 
M(n\!+\!r\!-\!m\!-\!1)$.  \\

An ordinary module is an $\N$-graded weak $V$-module 
$W=\oplus_{n=0}^{\infty} W(n)$ such that 
$L(0)$ acts on $W(n)$ semisimply and $\dim W(n)<\infty$ 
for all $n$. For a simple ordinary 
module $W=\oplus_{i=0}^{\infty}W(i)$, $L(0)$ acts on 
$W(0)$ as a scalar, which is called a 
{\bf conformal weight} of $W$. 
\end{dfn}

The main object in this paper is not an ordinary module, but 
an $\N$-graded weak module $W=\oplus_{m=0}^{\infty} W(m)$ 
with $\dim W(m)<\infty$. If $v\in V_m$, then $v^W(m\!-\!1)$ is 
a grade-preserving operator by 
(1') and we denote it by $o(v)$ and extend it linearly. 

\begin{dfn} $V$ is called {\bf $C_{2}$-cofinite} if the subspace 
$\left< v(\!-\!2)u: v,u\in V\right> $ has a finite codimension 
in $V$. 
\end{dfn}

Zhu has introduced the second vertex operator algebra 
$(V,Y[, ], \1, \tilde{\omega})$ associated to $V$ in 
Theorem 4.2.1 of \cite{Zh}. 

\begin{dfn}
The vertex operator $Y[v,z]=\sum_{n\in \Z}v[n]z^{-n-1}$ is defined 
for homogeneous $v$ via the equality
$$  Y[v, z]=Y(v, e^{z}\!-\!1)e^{z|v|} \in \End(V)[[z, z^{-1}]]  $$
and Virasoro element $\tilde{\omega}$ is defined to be 
$ \omega\!-\!\frac{c}{24}\1$, 
where $|v|$ denotes the weight of $v$. \\
\end{dfn}

Throughout this paper, we assume that 
$V=\oplus_{n=0}^{\infty}V_n$ 
is a $C_2$-cofinite VOA.

\subsection{$n$-th Zhu algebras}
Following \cite{FZ}, $V$ has a product 
$$v\ast u=\Res_x\frac{(1\!+\!x)^{|v|}}{x}Y(v, x)u  \eqno{(2.1)}$$
for $v \in V_{|v|}$ and $u \in V$, where $|v|$ denotes the weight 
of $v$.  Set 
$$O(V)=\left< \Res_x\frac{(1\!+\!x)^{|v|}}{x^2}Y(v, x)u 
\mid  v,u\in V\right> \eqno{(2.2)}$$ 
and $A(V)=V/O(V)$.  
Then it is known (Theorem 1.5.1 \cite{FZ}) that $A(V)$ is an 
associative algebra with a product $\ast$. We call it Zhu algebra. 
Zhu has also shown that $\omega+O(V)$ is in 
the center of $A(V)$. From now on, abusing the notation, 
we use the same notation 
$\omega$ for $\omega+O(V)$. 
The essential property of Zhu algebra is that a top 
module $W(0)$ of an 
$\N$-graded weak $A(V)$-module $W=\oplus_{n=0}^{\infty} W(n)$ is an 
$A(V)$-module and every $A(V)$-module is a top module of some 
$\N$-graded weak $V$-module. This concept was naturally extended 
to $n$-th graded piece of $\N$-graded weak modules 
by Dong, Li and Mason in \cite{DLiM2}. Set 
$$O_n(V)=
\left< \Res_x\frac{(x\!+\!1)^{|v|+n}}{x^{2+2n}}Y(v, x)u 
\mid v,u\in V\right>  
\eqno{(2.3)}$$ 
and $A_n(V)=V/O_n(V)$.  Like $A(V)$, $A_n(V)$ is an associative 
algebra with a product 
$$v\ast_nu=\sum_{m=0}^n\binom{-n}{m}
\Res_xY(v,x)u\frac{(x\!+\!1)^{|v|+n}}{x^{n+m+1}}$$ 
and has a property that an $n$-th (and less ) graded piece $W(n)$ 
of an $\N$-graded weak $V$-module $W=\oplus_{n=0}^{\infty}W(n)$ 
is an $A_n(V)$-module and every $A_n(V)$-module is an 
$n$-th (or less) graded piece of an $\N$-graded weak $V$-module, 
(see Theorem 4.2 in \cite{DLiM2}). $A_n(V)$ is called an 
$n$-th Zhu algebra, in particular, $0$-th Zhu 
algebra is the original Zhu algebra. It is easy to see that 
there is a natural homomorphism from $A_n(V)$ to $A_{n-1}(V)$. 
The product $\ast_n$ is characterized by the identity 
$o(v\ast_nu)=o(v)o(u)$ on $\oplus_{i=0}^nW(i)$ for every 
$\N$-graded weak $V$-module $W$ so that 
$A_n(V)$ is essentially the algebra of zero modes 
(grade-preserving operators) 
of fields. In particular, 
$\omega$ is a central element of $A_n(V)$ for any $n$. 
Viewing $A_n(V)$ as an algebra 
of zero modes, we will use the following notation:
$$  o(\al)=v(|v|\!-\!1\!+\!m)u(|u|\!-\!1\!-\!m) 
\quad \mbox{ in }\quad A_n(V)   $$ which implies that 
$o(\al)w=v(|v|\!-\!1\!+\!m)u(|u|\!-\!1\!-\!m)w$ 
for any $\N$-graded weak module $W$ and 
$w\in \oplus_{i=0}^nW(i).$ \\

We note that Garberdiel and Neitzke showed in \cite{GN} that 
if $V$ is a $C_2$-cofinite VOA of CFT type, then $A_n(V)$ is 
of finite dimension for 
any $n$. They didn't mention it directly. Later, Buhl proved 
it for irreducible modules in \cite{Bu} by using a spanning set 
of module. We would like to explain their results and 
prove the finiteness of dimension of $A_n(V)$ without 
assuming condition of being of CFT type. 
First, they showed that if $V=C_2(V)+\langle A\rangle $ for a 
set $A$ of $V$, then 
$V$ is spanned by elements of the form
$$  v^1(-N_1)\cdots v^r(-N_r)\1 $$
with $N_1\!>\!\cdots\!>\N_r\!>\!0$ and $v^i\in A$ 
by using a filtration, see Proposition 8 in \cite{GN}. 
We will first show the existence of such a 
spanning set for general VOAs without negative weights.  

\begin{lmm}
Let $V=\oplus_{m=0}^{\infty}V_m$ be a $C_2$-cofinite VOA and 
$A$ is a set of homogenous elements satisfying 
$V=C_2(V)+\langle A\rangle$. Let $W$ be a weak module generated 
from $w$. Then $W$ is spanned by 
the following elements  
$$  \{v^1(i_1)\cdots v^k(i_k)w  \mid v^i\in A, 
\quad i_1<\!\cdots\! <i_k \}.$$
In particular, if we set 
$$\CW(m)=\left< v^1(i_1)\cdots v^k(i_k)w \ \mid \ v^i\in V, \quad  
\deg(v^1(i_1)\cdots v^k(i_k))=m \right>,  $$
then $\CW(m)=0 \mbox{  for  } m\ll 0$, where 
$\deg(v(i)\cdots u(j))$ denotes the degree of $v(i)\cdots u(j)$ 
as an operator. 
\end{lmm}

\pr
Define a filtration on $W$ by  
$$\CW(n,m,r)=\left< v^1(i_1)\cdots v^k(i_k)w   
\left| \begin{array}{l}
v^i\in V, \quad  \sum_{i=1}^k\wt(v^i)\leq n, \cr 
\vspace{-4mm}\cr
\deg(v^1(i_1)\!\cdots\! v^k(i_k))=m, \quad k\leq r, \cr
 \end{array}\right.  \right> . $$  
Clearly, $\CW(n,m,r)\subseteq \CW(n\!+\!1,m,r)$ and 
$$W=\sum_{m\in \Z}\left(\cup_{n=0}^{\infty}
\cup_{r=0}^{\infty}\CW(n,m,r)\right).$$
We will prove that $\CW(n,m,r)$ is spanned by the desired 
elements contained in $\CW(n,m,r)$ for each $m$. 
Suppose false and let $(n,r)$ be a minimal counterexample 
with respect to lexicographical order.   
Let $U(n,m,r)$ be the 
subspace spanned by the desired elements contained in 
$\CW(n,m,r)$, 
then there is a nonzero element 
$u\in \CW(n,m,r)-U(n,m,r)$ and 
$\CW(n\!-\!1,m,r)+\CW(n,m,r\!-\!1)\subseteq U(n,m,r)$. 
We may assume $$u=v^1(i_1)\cdots v^r(i_r)w.$$
Since $v^i\in V\!=\!\langle A\rangle\!+\!C_2(V)$, there are $u^i\in A$ and 
$a^{ij},b^{ij}\in V$ such that $v^i=u^i+
\sum a^{ij}(\!-\!1)b^{ij}$. Since $C_2(V)$ and $\langle A\rangle$ are direct 
sums of homogeneous spaces, we may assume 
$\wt(v^i)=\wt(u^i)=\wt(a^{ij}(\!-\!2)b^{ij})
=\wt(a^{ij})\!+\!\wt(b^{ij})\!-\!1$. 
Using associativity 
$$(a(\!-\!2)b)(s)=\sum_{i=0}^{\infty} (\!-1\!)^i\binom{-2}{i}
\left\{ a(\!-\!2\!-\!i)b(s\!+\!i)-b(\!-\!2\!+\!s\!-\!i)a(i)
\right\},$$ 
we may assume 
$$  u=v^1(i_1)\cdots v^r(i_r)w \quad \mbox{ with } 
\quad v^i\in A.  \eqno{(2.4)}$$
We choose $u\in \CW(n,m,r)-U(n,m,r)$ so that $u$ has 
a form (2.4) and 
${\rm min}\{i_1,\cdots,i_k\}$ is minimal. The existence of 
minimal one follows from the next arguments. 
Since $v(i)u(j)=u(j)v(i)+[v(i),u(j)]=u(j)v(i)+
\sum \binom{i}{s}(v(s)u)(i\!+\!j\!-\!s)$ 
and $\wt(v(s)u)<\wt(v)+\wt(u)$ for $s\geq 0$, we may assume 
$$ u=v^1(i_1)\cdots v^r(i_r)w \quad \mbox{ with }\quad v^r\in A 
\mbox{  and  } i_1\leq \cdots \leq i_k. $$
Since $v^i\in A$ and $A$ is a finite set, there is an integer 
$N$ such that $v(m)w=0$ for $m>N$ and $v\in A$. 
Therefore we have $i_r>N$. In particular, 
the degree of $v^p(i_p)$ as an operator is bounded below. 
Since the total degree of $v^1(i_1)\cdots v^r(i_r)$ is $m$, 
the degree of $v^p(i_p)$ as an operator is bounded above  
and so $i_1$ is bounded below. 
By induction on $r$, we may assume $i_2<\cdots <i_r$. 
If $i_1<i_2$ or $i_1>i_2$, then $u\in \CW(n,m,r)$ by the 
minimality of $i_1$. The remaining case is $i_1=i_2$. 
Set $i=i_1=i_2$. The expansion of $(v^1(\!-\!1)v^2)_{2i+1}$ 
by associativity if 
$i\leq \!-\!1$ and that of $(v^2(\!-\!1)v^1)_{2i+1}$ 
by associativity if $i\geq 0$, 
contains a nonzero term $v^1_{i}v^2_i$ and the other terms are 
$v^1(i\!-\!j)v^2(i\!+\!j)$ or $v^2(i\!-\!j)v^1(i\!+\!j)$ 
with $j\not=0$. 
For example, if $i\leq \!-\!1$, then there are constants 
$\la_j$, $\mu_j$ such that 
$$v^1(i)v^2(i)=(v^1(\!-\!1)v^2)(2i\!+\!1)
+\sum_{j\not=0}( \la_jv^1(i\!-\!j)v^2(i\!+\!j)
+\mu_jv^2(i\!-\!j)v^1(i\!+\!j)).
\eqno{(2.5)}$$
We substitute the right-hand side for $v^1(i_1)v^2(i_1)$ in $u$. 
Since $v^1(i\!-\!j)v^2(i\!+\!j)v^3(i_2)\!\cdots\! v^r(i_r)w$ 
has a form (2.4) and 
one of $i\!-\!j$ and $i\!+\!j$ is less than $i$, we may assume 
$u=(v^1(\!-\!1)v^2)(2i\!+\!1)v^2(i_3)\cdots w$, which 
is in $\CW(n,m,r\!-\!1)$. Therefore we have a contradiction. 
Since $\dim \langle A\rangle<\infty$, 
there is an interger $N$ such that $i_k<N$ and so 
$\deg(v^1(i_1)\cdots v^k(i_k))
\geq \!-\!kN+\frac{k(k-1)}{2}>\!-\!N^2$ 
if $i_1<\cdots i_k<N$. 
Thus $\CW(m)=0$ if $m<\!-\!N^2$. 

This completes the proof. 
\prend

Using this spanning set, we have the following theorem, 
(see Lemma 3 and Theorem 11 in \cite{GN} if $V$ is of CFT type.)

\begin{thm}
If $V$ is $C_2$-cofinite, then $A_n(V)$ are all 
finite dimensional. 
\end{thm}

\pr
We fix $n$ and $v$ and $u$ denote homogeneous elements.  
Let $O_{(\infty^{2n+2})}$ be the subspace of $V$ spanned by 
elements of the form 
$$\left<  v(\!-\!2\!-\!N\!-\!2n|v|)u \ 
\mid \ v,u\in V, \ |v|\geq 1, 
\ N\geq 0 \right>.$$
We will show that $O_{(\infty^{2n+2})}$ is of finite codimension. 
Since $V$ is $C_2$-cofinite, we can choose a finite set $A$ of 
homogenous elements such that $V$ is spanned by 
$w_{1}(\!-\!N_1)\cdots w_{r}(\!-\!N_r)\1$  
with $N_1>\!\cdots\!>N_r$ 
and $w_{i}\in A$.  Let $t$ be the maximal weight of elements 
in $A$. If $N_1\geq (2n)t+2$, then $w_{i_1}(\!-\!N_1)\cdots 
w_{i_r}(\!-\!N_r)\1\in O_{(\infty^{2n+2})}$. 
This leaves 
us only finitely many choices for the $N_i$, which gives a finite 
spanning set for 
$V/O_{(\infty^{2n+2})}$. 

Let $O_{{\bf u}}$ be the subspace of $V$ spanned by 
elements of the form 
$$v\diamond_Mu=\dsp{\Res_z Y(v,z)u
\frac{(z\!+\!1)^{(n+1)|v|}}{z^{2n|v|+2+M}}}$$ 
with $|v|\geq 1$ and $M\geq 0$.  
Since $2n(|v|\!-\!1)\geq n(|v|\!-\!1)$ and 
$$v\diamond_Mu=\Res_z
Y(v,z)u\frac{(z\!+\!1)^{|v|+n+n(|v|-1)}}{z^{2+M+2n+2n(|v|-1)}},$$ 
 $v\diamond_Mu\in O_n(V)$ by Lemma 2.1 in \cite{DLiM2} and so 
it is sufficient to show \\
$\dim V/O_{{\bf u}}\!<\!\dim V/O_{(\infty^{2n+2})}$. We note 
$\dsp{v\diamond_M u
=\sum_{i=0}^{(n+1)|v|}\binom{(n+1)|v|}{i}
v(\!-\!2n|v|\!-\!2\!-\!M\!+\!i)u}$, 
that is, the weights of all terms are less than or equal to 
the weight of $v(\!-\!2n|v|\!-\!2\!-\!M)u$. 
Let $\{v_1,...,v_N\}$ be a set of representatives for 
$V/O_{(\infty^{2n+2})}$. 
Since $O_{(\infty^{2n+2})}$ is a direct sum of 
homogenous spaces, we may assume that 
$v_i$ are all homogenous elements. We claim that 
$\left<v_1,....,v_N\right> \!+\!O_{{\bf u}}=V$, which offers the 
desired conclusion. Suppose false and let 
$u\in V\!-\!\{\left< v_1,...,v_N\right> \!+\!O_{{\bf u}}\}$ 
be a homogeneous element with minimal weight among them. 
By the choice of $\{v_i\}$, there are 
$a_r, b_j\in \C$ and homogeneous elements $v^r, u^r\in V$ 
and $N_r\in \Z_+$ such that 
$$ u=\sum_j b_jv_j+\sum_r a_r v^r(\!-\!N_r\!-\!2n|v|\!-\!2)u^r.$$
We may assume that the weights of all elements in the above 
equation are the same. But then 
$$\hat{u}=u\!-\!\sum_j b_jv_j\!-\!\sum_r v^r\diamond_{N_r}u^r $$
is a linear combination of vectors whose weights are strictly 
smaller than that of $u$. 
By the minimality of $u$, $\hat{u}$ is contained 
in $\left< v_1,...,v_n\right> +O_{{\bf u}}$ 
and so is $u$. 

This completes the proof of Theorem 2.5. 
\prend

\subsection{Generalized Verma module}
In this paper, we will define a pseudo-trace function for a 
generalized Verma module $W^n_T$ constructed from an 
$A_n(V)$-module $T$. A generalized Verma module $\Verma(X)$ 
for an $A(V)$-module $X$ is 
introduced in \cite{Li1} as an extension of 
concept of a Verma module. 
It is a largest $\N$-graded weak $V$-module $W$ which has $X$ 
as a top module $W(0)$ and is 
generated from $X$. For an $A_n(V)$-module $T$, 
it is possible to consider a largest $\N$-graded weak 
$V$-module $W^{(n)}_T=\oplus_{i=0}^{\infty}W^{(n)}_T(i)$ 
which has $T$ as 
its $n$-th graded piece $W^{(n)}_T(n)$ 
and is generated from $T$. 
We should note that $W^{(n)}_T(0)$ 
might be zero if $T$ is also an $A_{n-1}(V)$-module. 
However, 
in this paper, we will only treat an $\N$-graded weak 
$V$-module $W^{(n)}_T$ 
constructed from an $A_n(V)$-module $T$ which 
satisfies the following 
condition: \\ 

\noindent
{\rm (2.6)}\quad For any $A_n(V)$-submodules 
$T^1\subseteq T^2$ of $T$ with 
an irreducible 
factor $T^2/T^1$, every $V$-module $W$ with $n$-th 
graded piece $T^2/T^1$ is irreducible. \\

More precisely, we will only consider an $\N$-graded weak 
$V$-module $W^{(n)}_T$ 
whose composition series has the same shape as does 
a composition series of $T=W^{(n)}_T(n)$. 
Therefore, a $V$-module $W$ with $T$ as an $n$-th 
graded piece is uniquely 
determined and so a generalized Verma module 
coincides with $L_n(T)$, which is defined 
in Theorem 4.2 of \cite{DLiM2} and is the minimal one in a sense. 
So we don't need a concept of generalized Verma module, 
but in order to emphasize that this 
module is naturally constructed from $A_n(V)$-module $T$, 
we will call it a generalized Verma module and denote it 
by $W^{(n)}_T$. As they showed in \cite{GN}, we obtain:

\begin{lmm} If $V$ is a $C_{2}$-cofinite VOA and 
$T$ is a finite dimensional $A_n(V)$-module, 
then $W_T(m)$ has a finite dimension for any $m=0,1,...$. 
\end{lmm}

\subsection{Elliptic functions}
We adopt the same notation form \cite{Zh}. 
The Eisenstein series $G_{2k}(\tau)$ $(k=1, 2, ...)$ are series 
$$G_{2k}(\tau)
=\sum_{(m, n)\not=(0, 0)}\frac{1}{(m\tau+n)^{2k}}  
\mbox{   for }k\geq 2 \qquad \mbox{  and} $$ 
$$G_2(\tau)=\frac{\pi^2}{3}
+ \sum_{m \in \Z-\{0\}}\sum_{n \in \Z}\frac{1}{(m\tau+n)^2} 
\mbox{   for }k=1.  $$
They have the $q$-expansions 
$$ G_{2k}(\tau)=2\xi(2k)+ \frac{2(2\pi\imath )^{2k}}{(2k\!-\!1)!}
\sum_{n=1}^{\infty}\frac{n^{2k\!-\!1}q^n}{1-q^n}, $$ 
where $ \xi(2k)=\sum_{n=1}^{\infty}\frac{1}{n^{2k}}$ and 
$q=e^{2\pi \imath \tau}$.
We make use of the following normalized Eisenstein series:
$$  E_k(\tau)=\frac{1}{(2\pi\imath )^k}G_k(\tau) 
\mbox{    for } k\geq 2. $$

\subsection{$C_2$-cofiniteness} 
Li showed in \cite{Li2} that $C_2$-cofinite VOA is finitely 
generated and regularity implies rationality 
and $C_2$-cofiniteness. Conversely, 
Abe, Buhl and Dong proved in \cite{ABD} that regularity comes from 
$C_2$-cofiniteness and rationality.   
The regularilty means that every weak module is a direct 
sum of simple ordinary modules. This implies two conditions, 
that is, every weak module is an $\N$-graded weak module and 
$V$ is rational. 
In this section, we will show that they are 
equivalent separately. Namely, $C_2$-cofinitenss 
implies that every weak module is an $\N$-graded weak module. 
We will prove the following theorem. 

\begin{thm}
Let $V=\oplus_{m=0}^{\infty}V_m$ be a vertex operator algebra. 
Then the following are equivalent. \\
{\rm (1)} \quad $V$ is $C_2$-cofinite. \\
{\rm (2)} \quad Every weak module is a direct sum of 
generalized eigenspaces of $L(0)$. \\
{\rm (3)} \quad Every weak module is an 
$\N$-graded weak module $W=\oplus_{n=0}^{\infty} W(n)$ such that $W(n)$ is a direct 
sum of generalized eigenspaces of $L(0)$. \\
{\rm (4)} \quad $V$ is finitely generated and every weak module 
is an $\N$-graded weak module. 
\end{thm}

\pr
Proof of (1) $\Rightarrow$ (3)  
(and so (1) $\Rightarrow$ (2) and (1)$ \Rightarrow$ (4) )\\
The proof is essentially the same as in \cite{ABD}. 
Let $W$ be a weak module. 
We note that $V$ has only finitely many irreducible modules. 
Let $\{r_1,...,r_k\}$ be the set of conformal weights and 
set $R=\cup_{i=1}^k( r_i\!+\!\Z_{\geq 0})$.
For $r\in \C$, $W_r=\{w\in W \mid (L(0)\!-\!r)^Nw=0 
\mbox{  for  }N\gg 0 \}$ 
denotes a generalized eigenspace of $L(0)$ with eigenvalue $r$. 
Clearly, $v(|v|\!-\!1\!+\!h)W_r\subseteq W_{r-h}$ for $v\in V$. 
We will show $W=\oplus_{r\in R}W_r$, which implies (3).  
It is easy to 
see that there is a unique maximal $\N$-graded weak submodule 
$U=\oplus_{n=0}^{\infty}U(n)$ 
such that $U(n)$ is a direct sum of generalized eigenspaces 
of $L(0)$. 
Clearly, $U=\oplus_{r\in R}U\cap W_r$. 
Suppose $W/U\not=0$. Since $W/U$ is also a weak module,
we may assume $W/U$ is generated from $w\not=0$. 
We use a similar argument as in \cite{ABD}. 
Applying Lemma 2.4 to $W/U$, there is an integer $m$ such that 
$\CW(m)\not=0$, but $\CW(t)=0$ for $t<m$. Then $\CW(m)$ is an 
$A(V)$-module. Since $A(V)$ is a finite dimensional algebra, 
$\CW(m)$ is a direct sum of generalized eigenspaces of $L(0)$, 
which contradicts the choice of $U$.  \\
Proof of (3) $\Rightarrow$ (2) is clear. \\
Proof of (2) $\Rightarrow$ (1). \\
As Li explained in \cite{Li2}, $V^{\ast}=\Hom(V,\C)$ contains 
a uniquely maximal weak submodule $W(V^{\ast})$ and if 
$f(C_2(V))=0$, then $f\in W(V^{\ast})$. 
By the assumption, $W(V^{\ast})$ is a direct sum of generalized 
eigenspaces of $L(0)$. 
Suppose $\dim V/C_2(V)=\infty$.  
Let $S=\{i \mid V_{i}\not=(C_2(V))_{i}\}$. For $i\in S$, 
choose $v_i\in V_i\!-\!(C_2(V))_i$ and 
a hyperspace $T_i$ of $V_i$ containing $(C_2(V))_i$ such that 
$V_i=\C v_i+T_i$. Set $T_i=(C_2(V))_i$ for $i\not\in S$. 
So $V$ is spanned by $T=\oplus T_i$ and $\{v_i\mid i\in S\}$.
Define $f\in V^{\ast}$ so that $f(T)=0$ and $f(v_i)=1$. 
Since $(L(0)f)(v_i)=f(L(0)v_i)=f(|v_i|v_i)=|v_i|$, 
$\{L(0)^if \mid i=0,1,... \}$ is a linearly independent, 
On the other hand, since $f\in W(V^{\ast})$, 
$f$ is a sum of finite elements in 
generalized eigenspaces of $L(0)$ and so 
$\left< L(0)^if \mid i=0,1,... \right>$ is of finite dimension,  
which is a contradiction. \\
Proof of (4) $\Rightarrow$ (1). \\
Since $V$ is finitely generated, $V/C_2(V)$ is a finitely 
generated abelian group (product is given by $u(\!-\!1)v$.) 
Therefore, there is a torsion-free element $u$, that is, 
$\{u(\!-\!1)^m\1\not\in \C_2(V)$ for any $m\in \N$. 
Define $f\in V^{\ast}$ such that $f(C_2(V))=0$ and 
$f(u(\!-\!1)^m\1)=1$ for any $m\in \N$. 
Then $v(2|v|\!-\!1)^mf\not=0$ for any $m\in \N$, which implies 
that $f$ does not belong to an $\N$-graded weak module. \\
This completes the proof of Theorem 2.7.

\section{Symmetric algebras and pseudo-trace maps}
In this section, we always consider a finite dimensional 
algebra over $\C$ with a unit $1$. 
Let $A$ be a ring and let $L(a)$ and $R(a)$ denote the left 
and right regular 
representations of $a\in A$ given by a basis $\{v_i \}$ of $A$. 
$A$ is called a Frobenius algebra 
if $L(a)$ and $R(a)$ are similar:  $L(a)=Q^{-1}R(a)Q$ for 
some matrix $Q$. In particular, $A$ is called a 
{\bf symmetric algebra} when the matrix $Q$ can be chosen as a 
symmetric matrix. It is also well known 
(cf. \cite{CR}) that this is equivalent to that 
$A$ has a symmetric linear map $\phi:A \to \C$ such that 
${\rm Rad}(\phi)\overset{{\rm def}}{=}
\{a\in A \mid \phi(AaA)=0 \}$ 
is zero and is also equivalent 
to that $A$ has a symmetric associative nondegenerated 
bilinear form 
$\langle \cdot,\cdot\rangle$, a relation is given by 
$\phi(ab)=\langle a,b\rangle$, where a symmetric map implies 
$\phi(ab)=\phi(ba)$ and an associative bilinear form means 
$\langle ab,c\rangle=\langle a,bc\rangle$ for $a,b,c\in A$, 
see \cite{CR}.
We denote a symmetric algebra $A$ with a symmetric linear 
function $\phi$ by $(A,\phi)$. \\

Originally, Nesbitt and Scott showed in \cite{NS} 
that $A$ is a symmetric 
algebra if and only if its basic algebra $P$ is symmetric. 
Later, Oshima gave a simpler proof of this 
equivalence, see \cite{O}. This is an idea to explain our 
strategy in this paper. Given an symmetric linear map $C_n$ of 
$n$-th Zhu algebra $A_n(V)$, we have a symmetric algebra 
$A=A_n(V)/\Rad(\phi)$ and 
its symmetric basic algebra $P$ with a symmetric linear 
function $\phi$. We start from $(P,\phi)$. We next 
construct a right $P$-module $W$ (a generalized Verma module) 
such that the basic algebra of $R=\End_P(W)$ is $P$.  
Nesbitt and Scott's result tells that $R$ has a symmetric 
linear map (we will call it pseudo-trace) and we will 
define a pseudo-trace function. 
In order to construct such right $P$-modules $W$, we will 
introduce a concept of modules interlocked with $\phi$ and 
the main purpose in this section is to construct a symmetric 
linear function $\tr^{\phi}_W$ of $R=\End_P(W)$ 
explicitly. \\

\noindent 
{\bf Setting}\\
Let $P$ be a basic symmetric algebra, that is, $P/J(P)$ is a 
direct sum of the field of complex numbers and $P$ has a symmetric 
linear function $\phi$ with $\Rad(\phi)=0$, 
where $J(P)$ denotes the Jacobson radical of $P$. Set \\ 

\vspace{-4mm}
(3.1) \quad $P/J(P)=\C \bar{e}_1\oplus 
\cdots\oplus \C \bar{e}_k$.\\

\vspace{-4mm}
(3.2) \quad $\{e_i\in P \ \mid \ i=1,...,k\}$ 
is the set of mutually orthogonal primitive idempotents \\ 

\vspace{-4mm}
\noindent
and $\bar{e}_i$ is the image of $e_i$ in $P/J(P)$. 
In particular, $1=e_1\!+\!\cdots\!+\!e_k$ is the identity of $P$. 
We also assume that: \\ 

\vspace{-4mm}
(3.3) \quad $P$ contains a central element $\omega$. \\ 

\vspace{-4mm}
\noindent
As we mentioned above, $P$ has a 
symmetric associative nondegenerated bilinear form 
$\langle\cdot,\cdot\rangle$ given by $\langle a,b\rangle
=\phi(ab)$. Since $\Rad(\phi)=0$, we obtain 
$\phi_{|M}\not=0$ for any minimal ideal $M$. On the other hand, 
there is a symmetric linear map $\pi$ 
of $P/J(P)$ such that $\pi(e_i)=\phi(e_i)$ for all $i$. 
Since such a map is a linear sum of ordinary trace maps, 
we may assume that we have started from $\phi-\pi$ and consider 
only the following case. \\ 

\vspace{-4mm}
(3.4) \quad $\phi(e_i)=0$ for all $i$.  \\ 

\vspace{-4mm}
\noindent
If $P$ is decomposable as a ring, 
say $P=P_1\oplus P_2$, then $(P,\phi)$ is a sum of 
$(P_1,\phi_{|P_1})$ and $(P_2,\phi_{|P_2})$. 
We hence assume that $P$ is indecomposable. In particular, 
there is $r\in \C$ and $\mu(r)\in \Z_+$ such that \\ 

\vspace{-4mm}
(3.5) \quad $(\omega\!-\!r)^{\mu(r)}P=0$. \\

First of all, we have: 

\begin{lmm}
$\soc(P)\cong \C e_1\oplus \cdots\oplus \C e_k$ 
as $P\times P$-modules 
and $\soc(P)\subseteq J(P)=\soc(P)^{\perp}$, 
where $\soc(P)$ denotes the socle of $P$, that is, the 
sum of all minimal left ideals.  
\end{lmm}

\pr
Let $M$ be a minimal ideal. If $e_iMe_j=M$ for $i\not=j$, then 
$\phi(m)=\phi(e_im\!-\!me_i)=0$ for all $m\in M$, which 
contradicts $\Rad(\phi)=0$. 
Hence $M$ satisfies 
$e_iMe_i=M$ for some $i$.  For $i=1,...,k$, 
$M_i=(\C e_1+\cdots+\C e_{i-1}+\C e_{i+1}+\C e_k+J(P))^{\perp}$ 
is a minimal ideal with $e_iM_ie_i=M_i$.
If $P$ has two minimal ideals $M$ and $N$ which are isomorphic 
together, then 
$L=\{(a,b)\in M\oplus N  \mid \phi(a)+\phi(b)=0 \}$ is a nonzero 
ideal of $P$ with $\phi(L)=0$, which is a contradiction. 
Hence $\soc(P)\cong \C e_1\oplus\cdots \oplus \C e_k=J(P)^{\perp}$. 
Since $P$ is indecomposable, $\soc(P)\subseteq J(P)$. 
\prend

Let $\{f_i\in \soc(P) \mid i=1,...,k\}$ be the dual basis of 
$\{e_i \mid i=1,...,k\}$. 
We note $\phi(f_i)=\langle 1,f_i\rangle
=\langle e_i,f_i\rangle=1$. 
Set $d_{ij}=\dim_{\C} e_iJ(P)e_j/e_i\soc(P)e_j$.  
Then we can find the following basis of $P$. 

\begin{lmm} 
$P$ has a basis 
$$\Omega=\{\rho^{ii}_0, \ \rho^{ii}_{d_{ii}+1}, 
\ \rho^{i,j}_{s_{ij}} \  \mid 
\ i,j=1,...,k, \ s_{ij}=1,...,d_{ij}\}$$ 
satisfying \\
{\rm (1)}\quad $\rho^{ii}_0=e_i$, $\rho^{ii}_{d_{ii}+1}=f_i$, \\
{\rm (2)}\quad $e_i\rho^{ij}_se_j=\rho^{ij}_s\quad $  
for all $i,j,s$, \\
{\rm (3)}\quad $\langle \rho^{ij}_s, \rho^{ab}_{d_{ab}-t}\rangle
=\delta_{ib}\delta_{ja}\delta_{s,t}\quad $ and \\
{\rm (4)}\quad $\rho^{ij}_s\rho^{ji}_{d_{ji}-s}=f_i$. 
\end{lmm}

\pr
We first choose $\rho^{ii}_0=e_i$ and $\rho^{ii}_{d_{ii}+1}=f_i$, 
which satisfy (2), (3) and (4). As we showed, 
$J(P)^{\perp}=\soc(P)\subseteq J(P)$ and so $J(P)/\soc(P)$ has 
a nondegenerated symmetric bilinear form and so do 
$e_iJ(P)e_i/e_i\soc(P)e_i$ and 
$e_iJ(P)e_j/e_i\soc(P)e_j\!+\!e_jJ(P)e_i/e_j\soc(P)e_i$. 
Since $J(P)$ is nilpotent, there is a sequence of ideals 
$\FA_i$ of $e_iPe_i$:
$$  e_iPe_i=\FA_0\supsetneq J(P)\cap e_iPe_i
=\FA_1\supsetneq \cdots \supsetneq \FA_{d_{ii}+1}
=\soc(P)\cap e_iPe_i 
\supsetneq \FA_{d_{ii}+2}=0  $$
such that $\FA_s/\FA_{s+1}$ is a simple $P$-module and 
$\FA_s\FA_{d-s+1}\subseteq e_i\soc(P)e_i$. 
In particular, we may choose $\FA_s$ so that 
$\FA_s^{\perp}=\FA_{d_{ii}+2-s}$ for all $s$. 
Choosing a base 
$\rho_s^{ii}$ in $\FA_s\!-\!\FA_{s+1}$ and 
its dual base $\rho^{ii}_{d_{ii}-s+1}=(\rho^{ij})^{\ast}$ in 
$\FA_{s+1}^{\perp}$ for $s\leq \hf(d_{ii})$ so that 
$\langle \rho^{ii}_s,\rho^{ii}_{t}\rangle=\delta_{s+t,d_{ii}+1}$ 
inductively, 
we have a desired set for $e_iPe_i$. For $i\not=j$, 
we first note $d_{ij}=d_{ji}$. 
Then there are sequences of ideals $\FB_s$ and $\FC_s$ of 
$e_iPe_j$ and $e_jPe_i$: 
$$ \begin{array}{l}
 e_iPe_j=\FB_1\supsetneq \FB_2\supsetneq \cdots \supsetneq 
\FB_{d_{ij}} \supsetneq 0  \cr
 e_jPe_i=\FC_1\supsetneq \FC_2\supsetneq \cdots \supsetneq 
\FC_{d_{ij}} 
\supsetneq 0  
\end{array}$$
such that $\FB_s\FC_{d_{ij}-s+1}=\C f_i$ and 
$\FC_{d_{ij}-s+1}\FB_s=\C f_j$. 
So we can take $\rho^{ij}_s\in \FB_s\!-\!\FB_{s+1}$ and choose 
suitable $\rho^{ij}_{d_{ij}-s+1}\in 
\FC_{d_{ij}-s+1}\!-\!\FC_{d_{ij}-s}$ so that 
they satisfy (3) and (4). It is easy to see that the set of 
all $\rho^{ij}_s$ satisfies the desired conditions. 
\prend

For $\rho\in \Omega$, $\rho^{\ast}$ denotes its dual, 
that is, a unique 
element $\rho^{\ast}\in \Omega$ such that 
$\langle \rho,\rho^{\ast}\rangle=1$. 

\begin{lmm}
Set $P^0=(\C e_1+\cdots+\C e_k)^{\perp}$. Then 
for $\rho, \mu\in \Omega$, if $\mu\not=\rho^{\ast}$, 
then $\rho\mu\in P^0$. 
\end{lmm}

\pr 
Assume $e_j\rho=\rho$. Then we have 
$\langle e_i, \rho\mu\rangle=\langle e_i\rho,\mu\rangle
=\langle \delta_{ij}\rho, \mu\rangle
=\delta_{ij}\langle \rho,\mu\rangle=0$ for all $i$.  
\prend

We assume that $P$ has a basis satisfying 
the conditions in Lemma 3.2. As we explain, we would like to 
call a right $P$-module $W$ ``interlocked with $\phi$'' if a basic 
algebra of $\End_P(W)$ is $P$. However, we will need a symmetric 
linear map $\tr^{\phi}$ explicitly and so we will just 
introduce a sufficient condition for that.  

\begin{dfn}
Let $(P,\phi)$ be a basic symmetric algebra with a symmetric 
linear map $\phi$ satisfying (3.4) and $W$ is a finite 
dimensional right $P$-module. We will 
call that $W$ is {\bf interlocked with} $\phi$ if 
$\Ker(f_i)\overset{{\rm def}}{=}\{w\in W \mid  wf_i=0\}$ is equal 
to $\sum_{\rho\in \Omega-\{e_i\}} W\rho$ for each $i=1,...,k$, 
where $f_i=\rho^{ii}_{d_{ii}+1}$ is a dual base of $e_i$. 
\end{dfn}

Set $R=\End_P(W)$. Then 
$W\al$ is an $R$-module for each 
$\al\in P$ and the condition above 
implies that 
$Wf_i\cong W/\{w\in W  \mid wf_i=0\}
=W/\sum_{\al\in \Omega-\{e_i\}}W\al
\cong W/(\sum_{i\not=j} We_j+WJ(P))$ 
as $R$-modules. Namely, the definition says that the 
socle part $\oplus_{i=1}^kWf_i$ and 
the semisimple part 
$W/WJ(P)$ are isomorphic together. 
Not only the top and the bottom, but we also have: 

\begin{lmm}
$We_i/WJ(P)e_i \to We_i\rho e_j/WJ(P)e_i\rho e_j$ is also an 
$R$-isomorphism for any $e_i\rho e_j\in \Omega$. 
\end{lmm}

\pr
Since the composition map 
$$(e_i\rho e_j)^{\ast}(e_i\rho e_j): We_i/\sum_{\mu<e_i} W\mu 
\overset{e_i\rho e_j}{\longrightarrow }
We_i\rho e_j/\sum_{\mu<e_i}W\mu e_i\rho e_j 
\overset{(e_i\rho e_j)^{\ast}}{\longrightarrow}  
We_i\rho e_je_j\rho^{\ast}e_i=W f_i$$  
is an $R$-isomorphism, so is 
$$e_i\rho e_j: We_i/\sum_{\mu<e_i} W\mu \to 
We_i\rho e_j/\sum_{\mu<e_i}W\mu e_i\rho e_j$$ 
for any $e_i\rho e_j\in \Omega$, where $\mu<e_i$ implies 
$W\mu \subsetneq We_i$.  
\prend

Therefore, it is easy to see that 
$W$ is interlocked with $\phi$ if and only if there are 
vector spaces $T_p$ ($\cong Wf_p$) for each $p$ such that 
$$W\cong \oplus_{e_p\rho e_s\in \Omega} T_p\otimes e_p\rho e_s$$ 
as right $P$-modules. 
Now we will define a symmetric linear map $\tr^{\phi}_W$ 
of $R=\End_P(W)$. 
Let $\{v^{e_p}_i \mid i=1,...,\dim T_p \}$ be a basis of $T_p$. 

\begin{dfn}
Let $(P,\phi)$ be a basic symmetric algebra and $W$ a right 
$P$-module interlocked with $\phi$. 
For $\al\in R=\End_P(W)$, we have a 
$(\sum_{p=1}^k \dim T_p)\times
(\sum_{e_p\rho e_s\in \Omega}\dim T_p)$-matrix  
$\left(\al^{e_p\rho e_s}_{ji}\right)$ such that 
$$\al(v_j^{e_s}\otimes e_s)
=\sum_{e_p\rho e_s\in \Omega}
\left(\sum_{i=1}^{\dim T_p} \al_{ji}^{e_p\rho e_s}v_i^{e_p}
\otimes e_p\rho e_s \right)$$
for $v^{e_s}_j\otimes e_s$.  
Then we define a pseudo-trace map $\tr^{\phi}_W(\al)$ by the 
sum of traces of 
matrices $(\al_{ji}^{f_s})_{ji}$ for $s=1,...,k$. 
Namely, we define 
$$\tr^{\phi}_W(\al)
=\sum_{s=1}^k\tr (\al_{ji}^{f_s})
=\sum_{s=1}^k\sum_{i=1}^{\dim T_s}\al^{f_s}_{ii}.$$  
\end{dfn}

Namely, we take the trace of $\Hom(W/W J(P) \to W \soc(P))$ as 
an example in the introduction. 

\begin{prn}
$\tr^{\phi}_W$ is a symmetric linear map. 
\end{prn}

\pr
It is easy to see that $\tr^{\phi}_W$ does not depend on the 
choice of bases of $T_p$'s. 
Let $\al,\be\in R=\End_P(W)$. 
Then there are $\al_{ji}^{e_p\rho e_s}, 
\be_{ji}^{e_p\rho e_s}\in \C$ such that 
$$\begin{array}{l}
\dsp{\al(v_j^{e_s}\otimes e_s)=
\sum_{e_p\rho e_s\in \Omega}
\left(\sum_i \al_{ji}^{e_p\rho e_s}v_i^{_p}
\otimes e_p\rho e_s \right) } \qquad \mbox{ and }\cr
\dsp{\be(v_i^{e_p}\otimes e_p)=\sum_{e_t\mu e_p\in \Omega}
\left(\sum_h \be_{ih}^{e_t\mu e_p}v_h^{e_t}
\otimes e_t\mu e_p \right)}
\end{array}$$
By direct calculation, we obtain 
$$\begin{array}{rl}
\be\al(v_j^{e_s}\otimes e_s)=&
\be\left(\sum_{e_p\rho e_s\in \Omega}
\left(\sum_i \al_{ji}^{e_p\rho e_s}
v_i^{e_p}
\otimes e_p\rho e_s \right)\right) \cr
=&\sum_{e_p\rho e_s\in \Omega}
\left(\sum_i \al_{ji}^{e_p\rho e_s}\be(v_i^{e_p}
\otimes e_p)e_p\rho e_s \right) \qquad 
\mbox{ since } \be\in \End_P(W), \cr
=&\sum_{e_p\rho e_s\in \Omega}\left(\sum_i \al_{ji}^{e_p\rho e_s}
\left(
\sum_{t,\mu}\left(\sum_h \be_{ih}^{e_t\mu e_p}v_h^{e_t}
\otimes e_t\mu e_p \right)\right)e_p\rho e_s\right) \cr
=&\sum_{e_p\rho e_s\in \Omega}\sum_i\sum_{t,\mu}\sum_h 
\al_{ji}^{e_p\rho e_s}
\be_{ih}^{e_t\mu e_p} v_h^{e_t}\otimes 
e_t\mu e_pe_p\rho e_s.
\end{array}$$
By Lemma 3.3, $e_t\mu e_pe_p \rho e_s\in P^0
=\left< \Omega\!-\!\{f_1,...,f_k\}\right>$ except for 
$e_t\mu e_p=(e_p\rho e_s)^{\ast}$ and so we have  
$$\tr^{\phi}_W(\be\al)=\sum_{s=1}^k\tr (\al_{ji}^{e_p\rho e_s})
(\be_{ih}^{e_s\rho^{\ast}e_p}).$$
Hence we obtain 
$$
\tr^{\phi}_W(\be\al)=\sum_{s=1}^k\tr (\al_{ji}^{e_p\rho e_s})
(\be_{ih}^{e_s\rho^{\ast}e_p})=\sum_{s=1}^k
\tr (\be_{ih}^{e_s\rho^{\ast}e_p})
(\al_{ji}^{e_p\rho e_s})=\tr^{\phi}_W(\al\be),$$
as desired.
\prend

We will call $\tr^{\phi}_W$ a pseudo-trace map. We should 
note that we will 
treat an ordinary trace map as one of pseudo-trace map although 
$\phi$ does not satisfy (3.4). For if a basic 
symmetric algebra $(P,\phi)$ is $(\C, 1)$, 
then $\tr^1_W$ coincides with the ordinary trace map. 
We next investigate the properties of pseudo-trace maps. 
Let $\omega$ be a central element of $P$ such that 
$(\omega\!-\!r)^sP=0$ and $(\omega\!-\!r)^{s-1}P\not=0$ for 
$r\in \C$. Set $\FN=\{a\in P  \mid (\omega\!-\!r)a=0\}$. 
It is easy to see that if we define $\phi'$ by 
$$\phi'(a)=\phi((\omega\!-\!r)a),$$  
then $\phi'$ is also a symmetric linear function of $P/\FN$. 
We denote it by $(\omega\!-\!r)\phi$ and the right action of 
$\omega\!-\!r$ by $(\omega\!-\!r)_P$. 

\begin{prn}
If $W$ is a right $P$-module interlocked with $\phi$, then 
$W/W\FN$ is a right $P/\FN$-module interlocked with 
$(\omega\!-\!r)\phi$ and 
$$\tr^{\phi}_W \left( g(\omega\!-\!r)_P\right)
=\tr^{(\omega-r)\phi}_{W/W\FN} g  \eqno{(3.6)}$$
for $g\in End_P(W)$. (We also view $g\in \End_{P/\FN}(W/W\FN)$.)
\end{prn}

\pr
Set $R=\End_P(W)$. 
The first assertion is clear. Since $\omega\!-\!r\in Z(P)$, we 
obtain $\omega\!-\!r\in Z(R)$, where $Z(A)$ denotes the 
center of $A$.  Let $\FD$ be an idela of $R$ such that 
$\FD/R\FN=\soc(R/R\FN)$. Then 
$(\omega\!-\!r)\FD\subseteq \soc(R)\cap \Image(\omega\!-\!r)$ and 
$(\omega\!-\!r)\FD\cong \FD/\FN$. On the other hand, 
$\FE=\{a\in R \mid  (\omega\!-\!r)a\in \soc(R)\}$ satisfies 
$\FE/R\FN\subseteq 
\soc(R/R\FN)$. Hence $(\omega\!-\!r)$ is an isomorphism from 
$\soc(R/R\FN)$ to $\soc(R)\cap \Image(\omega\!-\!r)$. From the 
definition of pseudo-trace maps, 
$\tr^{\phi}_W$ is given by the traces of $\Hom_P(R/J(R),\soc(R))$ 
and so $\tr^{\phi}_W(\omega\!-\!r)$ is given 
by the traces of $\Hom_P(R/J(R), (\omega\!-\!r)^{-1}\soc(R)/R\FN)$, 
which equal the traces of $\Hom_P(R/J(R),\soc(R/\FN))$. 
We hence have 
$\tr^{\phi}_W \left( g(\omega\!-\!r)_P\right)
=\tr^{(\omega-r)\phi}_{W/W\FN} g$, 
as desired. 
\prend

Let $A$ be a finite dimensional symmetric algebra with a symmetric 
linear map $\phi$. Let $A/J(A)=A_1\oplus\cdots\oplus A_k$ 
be the decomposition into the direct sum of simple components. 
Let $\{e_i \mid i=1,...,k\}$ be a set of orthogonal 
idempotents such that 
$\bar{e}_i=e_i+J(A)$ is a primitive idempotent 
of $A_i$ for each $i$. Set $e=\sum_{i=1}^k e_i$. Then $eAe$ is 
a basic algebra and 
$eAe$ has a symmetric map $\phi$ (we use the same notation 
$\phi$) with zero radical. 
Viewing $Ae$ as a 
right $eAe$-module, it is easy to check that $Ae$ is interlocked 
with $\phi$ and so we can define a pseudo-trace map 
$\tr^{\phi}_{Ae}$ of 
$A\subseteq \End_{eAe}(Ae)$ on $Ae$. From the 
definition of $\tr^{\phi}_{Ae}$ and the previous arguments, 
it is easy to see the 
following.

\begin{lmm}
$\tr^{\phi}_{Ae}(a)=\phi(a)$ for all $a\in \soc(A)$.
\end{lmm}

So we have the following theorem, which we need later.

\begin{thm} Let $A$ be a finite dimensional associative algebra 
over ${\C}$ and $\phi$ a linear function of $A$ satisfying 
$\phi(ab)=\phi(ba)$ for every $a,b \in A$.  
Let $\omega$ be in the center of $A$ and assume 
that $\phi((\omega\!-\!r)^{\mu(r)}a)=0$ 
for every $a\in A$. Then there are linear 
symmetric functions $\phi_i$ of $A$  $(i=1,...,s)$ and 
basic symmetric subalgebras $P_i$ of factor rings 
$A/\FN_i$ with symmetric 
linear functions $\phi_i$ and $A\times P_i$-modules 
$M^i=(A/\FN_i)e$ 
satisfying $(\omega\!-\!r)^{\mu(r)}M^i=0$ such that 
$$  \phi(b)=\sum_{i=1}^s \tr^{\phi_i}_{M^i}(b) $$
for every $b \in A$, where $\FN_i=\Rad(\phi_i)$. 
\end{thm}

\pr 
We will use induction on $\dim A$. 
If $\phi$ has a nonzero radical $M$, then 
$\phi:A/M \to \C$ satisfies the 
same condition, but $\dim A/M< \dim M$. So we may assume that 
$\Rad(\phi)$ is zero and $A$ is a symmetric 
algebra with a symmetric linear map $\phi$. We may also 
assume that $A$ is indecomposable. 
By the previous lemma, there is an $A\times eAe$-module $Ae$ and 
a symmetric linear map $\psi$ such that 
$\tr^{\psi}_{Ae}(a)=\phi(a)$ for all $a\in \soc(A)$. 
Set $\phi'=\phi\!-\!\tr^{\psi}_{Ae}$. Then 
$\phi'$ satisfies the same condition and $\Rad(\phi')\supseteq 
\soc(A)$. By induction on $\dim A$, $\phi'$ is a sum of 
pseudo-trace maps and so is 
$\phi$. 
\prend

\section{Pseudo-trace functions}
Let $\{W^1,...,W^k\}$ be the set of all irreducible $V$-modules 
and set $W^i=\oplus_{m=0}^{\infty}W^i(m)$ with $W^i(0)\!\not=\!0$. 
It may happen that $W^i(h)=0$ for some $h>0$ and 
$W^i(h\!-\!1)\not=0$. Then, since 
$sl_2(\C)\cong \left< L(\!-\!1),2L(0),L(1)\right>$ acts on 
$\oplus_{m=0}^{h-1}W^i(m)$, 
the conformal weight is $1\!-\!h$. Therefore, 
there is an integer ${\it l}$ 
such that $W^i(m)\not=0$ for any $m>{\it l}$ and $i$.  

In this section, we assume: \\

\noindent
{\rm (4.1)}\quad $A_n(V)$ has a symmetric linear map $\phi$ 
satisfying $\phi((\omega\!-\!r)^s\ast v)=0$ for some 
$r\in \C$ and $s\in \N$ and 
the real part ${\rm Re}(r)$ of $r$ is greater than the real 
parts of conformal weights of all irreducible 
modules by ${\it l}$. \\

Set $A=A_n(V)/\Rad(\phi)$, which 
is a symmetric algebra with a 
symmetric linear map $\phi$. 
Let $A/J(A)=A_1\oplus\cdots \oplus A_k$ be the 
decomposition into the direct sum of simple components $A_i$. 
Let $\{e_i \mid i=1,...,k\}$ be a set of mutually 
orthogonal primitive idempotents 
of $A$ such that $\bar{e}_i=e_i+J(A)$ is a 
primitive idempotent of $A_i$. Set $e=e_1+\cdots+e_k$. 
Then $P=eAe$ is a basic symmetric algebra with 
a symmetric linear map 
$\phi$. 
As we showed in the previous section, 
$T=Ae$ is interlocked with $\phi$. 

Let $W^{(n)}_T$ denote a generalized Verma $V$-module 
generated from 
$A_n(V)$-module $T$. We assert that 

\begin{center}
$W^{(n)}_T$ is interlocked with $\phi$. 
\end{center}

We first note that $W^{(n)}_T(m)$ has a finite dimension 
for every $m$ 
and $L(0)\!-\!r$ acts on $W^{(n)}_T(n)\cong T$ as a nilpotent 
operator. By the definition of generalized Verma module 
$W^{(n)}_T$, $W^{(n)}_T$ is a right $P$-module by the action 
$$(\sum v^1(i_1)\cdots v^s(i_s)x )g
=\sum v^1(i_1)\cdots v^s(i_s)(xg)$$
for $v^i\in V$, $x\in T$ and $g\in P$. 
Since the real part ${\rm Re}(r)$ of $r$ is greater 
than the real parts of conformal weights of 
any modules by ${\it l}$, 
every nonzero submodule of $W^{(n)}_T$ has a nonzero intersection 
with $W^{(n)}_T(n)$ and so all irreducible factors $W/U$ of 
composition series of 
$W^{(n)}_T$ with 
isomorphic $n$-th graded pieces as $A_n(V)$-modules are 
isomorphic together as $V$-modules. 
In particular, the semisimple part $W^{(n)}_T/W^{(n)}_TJ(P)$ 
and the socle part 
$W^{(n)}_T\soc(P)$ are isomorphic together as $V\times P$-modules. 
Thus, $W^{(n)}_T$ is interlocked with $(P,\phi)$. \\

We note that $W^{(n)}_T(0)$ may be zero if $T$ is also an 
$A_{n-1}(V)$-module. We should also note that we have defined 
a pseudo-trace map for a finite dimensional vector space $W$ and 
so we have to say that for each $N$, $\oplus_{h=0}^NW^{(n)}_T(h)$ 
is interlocked with $\phi$. However, by the definition of 
pseudo-trace map, $\tr^{\phi}_W$ on $W=\oplus_{h=0}^NW^{(n)}_T(h)$ 
does not depend on the choice of $N$ and so it is uniquely defined 
on any $W_T^{(n)}(h)$ as does an 
ordinary trace map.  
We also note that 
$$\tr^{\phi}_{W^{(n)}_T(n)}(a)=\tr^{\phi}_{T}(a)=\phi(a) $$
for $a\in \soc(A)$ by definition.

\subsection{Logarithmic modules}
In this subsection, $\omega$ denotes Virasoro element and 
we fix $r\in \C$ and $s\in \N$. Assume $(\omega\!-\!n\!-\!r)^sT=0$ 
and $(\omega\!-\!n\!-\!r)T\not=0$, thus 
$L(0)=o(\omega)$ does not act on the graded piece $W^{(n)}_T(m)$ 
semisimply. 
However, since $(L(0)\!-\!m\!-\!r)^sW^{(n)}_T(m)=0$, we are 
able to understand 
$$e^{2\pi iL(0)\tau}=e^{2\pi i(m+r)\tau}
\left(\sum_{j=0}^{s-1}\frac{1}{j!}
(2\pi i\tau (L(0)\!-\!m\!-\!r))^j\right) 
\mbox{  on } W^{(n)}_T(m) $$ 
and define $q^{L(0)}$ on $W^{(n)}_T(m)$ by 
$$q^{m+r}(\sum_{j=0}^{s-1}\frac{1}{j!}
\left(2\pi i\tau (L(0)\!-\!m\!-\!r))^j\right). $$
Such a module is called a logarithmic module. 
We note that the left action $L(0)\!-\!m\!-\!r$ on $W^{(n)}_T(m)$ 
is equal to 
the right action of $\omega\!-\!n\!-\!r\in P$ on $W(m)$ and 
we denote it by $(\omega\!-\!n\!-\!r)_P$. 
Set $\FN_i=\{a\in P \mid (\omega\!-\!n\!-\!r)^ia=0\}$ and 
let $L^s(0)$ be a 
degree operator which acts on 
$W(m)$ as $m\!+\!r$, that is, the semisimple 
part of $L(0)$ and $L(0)\!-\!L^s(0)$ is nilpotent.  
By Proposition 3.8, we have: 

\begin{lmm}
$$\tr^{\phi}_W (L(0)\!-\!L^s(0))^ig=\tr^{\phi}_W 
g(\omega\!-\!n\!-\!r)^i_P
=\tr^{(\omega-n-r)^i\phi}_{W/W\FN_i}g. $$
\end{lmm}

Set 
$$ S^W(v,\tau)=\tr^{\phi}_W o(v)q^{L(0)-c/24}, 
\qquad (q=e^{2\pi i\tau})  $$
for a generalized Verma module $W=W^{(n)}_T$ 
interlocked with $\phi$. 
Although we are studying a general (nonsemisimple) 
operator $L(0)$ and 
a pseudo-trace function, 
they satisfy the following properties as do the action of $L(0)$ 
on modules and trace map: \\
(1)\quad $\tr^{\phi}_W$ is a symmetric function, \\
(2)\quad $[L(0), v(m)]=(|v|\!-\!m\!-\!1)v(m)$ and \\
(3)\quad $\tr^{\phi}_W o(\omega)o(v)q^{L(0)}
=\frac{1}{2\pi i}\frac{d}{d\tau}\tr^{\phi}_Wo(v)q^{L(0)}$, \\ 
which are the properties that Zhu used in the proof 
for ordinary trace functions. Therefore we have the 
following results 
by exactly the same arguments as in \cite{Zh}. 
 
\begin{prn}  For any $v,u \in V$ and a generalized Verma module 
$W$ interlocked with $(P,\phi)$, 
we have 
$$
\tr^{\phi}_W v(|v|\!-\!1\!-\!k)u(|u|\!-\!1\!+\!k)q^{L(0)}
=\tr^{\phi}_W
\frac{-q^k}{1-q^k}\sum_{s=0}^{\infty}\binom{|v|\!-\!1\!-\!k}{s}
o(v(s)u)q^{L(0)}, $$
$$\tr^{\phi}_W u(|u|\!-\!1\!-\!k)v(|v|\!-\!1\!+\!k)q^{L(0)}
=\tr^{\phi}_W
\frac{q^k}{1-q^k}\sum_{s=0}^{\infty}\binom{|v|\!-\!1\!+\!k}{s}
o(v(s)u)q^{L(0)},$$ 
$$ \tr^{\phi}_{W} o(v[0]u)q^{L(0)}=0,    $$
$$\tr^{\phi}_{W} o(v)o(u)q^{L(0)}
=\tr^{\phi}_{W} o(v[-1]u)q^{L(0)}-\sum_{k=1}^{\infty}E_{2k}
(\tau)\tr^{\phi}_{W} o(v[2k\!-\!1]u)q^{L(0)},$$
$$ \tr^{\phi}_W o(v[\!-\!2]u)q^{L(0)}+ \sum_{k=2}^{\infty}
(2k\!-\!1)E_{2k}(\tau){\tr^{\phi}}|_Wo(v[2k\!-\!2]u)q^{L(0)}=0, 
\mbox{  and } $$ 
$$\tr^{\phi}_{W} o(L[\!-\!2]u)q^{L(0)-c/24}
-\sum_{k=1}^{\infty}E_{2k}(\tau)\tr^{\phi}_{W} 
o(L[2k\!-\!2]u)q^{L(0)-c/24}
=\frac{1}{2\pi i}\frac{d}{d\tau}(\tr^{\phi}_{W} 
o(u)q^{L(0)-c/24}). $$
\end{prn}

\subsection{One-point functions}

By the investigation of trace functions in \cite{Zh}, Zhu showed 
$$\begin{array}{l}
\dsp{S^M(v[\!-\!2]u
+\sum_{k=1}^{\infty}(2k\!-\!1)E_{2k}(\tau)v[2k\!-\!2]u,\tau)=0 }\cr
\dsp{S^M(L[\!-\!2]u
+\sum_{k=1}^{\infty}E_{2k}(\tau)L[2k\!-\!2]u,\tau)
=\frac{1}{2\pi i}\frac{d}{d\tau}S^W(u,\tau)}
\end{array}$$
and 
$$S^M(v[0]u,\tau)=0 \eqno{(4.2)} $$
for a $V$-module $M$.  
Since $[o(v),o(u)]=o(v[0]u)$, (4.2) implies that $S_W(\ast,\tau)$ 
is a symmetric linear function on 
$\left< o(v) \mid v\in V\right>$ in a sense. 

Consider $V[E_4(q),E_6(q)]\subseteq V[[q]]$. 
$O_q(V)$ is the submodule of $V[E_4(q),E_6(q)]$ generated by 
elements of the type 
$$v[0]u$$ 
and 
$$  v[\!-\!2]u+ \sum_{k=2}^{\infty} (2k\!-\!1)E_{2k}(\tau)
\otimes v[2k\!-\!2]u    \quad v,u \in V. $$

We first prove the following lemma. 

\begin{lmm}
For $\al\in V$ and a fixed integer $n>{\it l}$, there are $v^i$ 
and $u^i \in V$ $(i=1,...,p)$ such that 
$$o(\al)=\sum_{i=1}^p v^i(|v^i|\!-\!1\!+\!n)u^i(|u^i|\!-\!1\!-\!n) 
\mbox{  in  }A(V),  $$
where the statement 
``$o(\be)=v^i(|v^i|\!-\!1\!+\!n)u^i(|u^i|\!-\!1\!-\!n)$ in $A(V)$'' 
implies $o(\be)u=v^i(|v^i|\!-\!1\!+\!n)u^i(|u^i|\!-\!1\!-\!n)u$ for 
any $\N$-graded weak $V$-module $W$ and $u\in W(0)$. 
\end{lmm}

\pr
We first note the following fact, which is a natural 
consequence of associativity.  
For any $v, u \in V$ and $r, s\in \Z$ and $m\in \N$, 
there is an element $\be$ of $V$ which is a 
linear combination $\be=\sum \la_iv(i)u+\sum \mu_i u(i)v$ 
of $v(i)u$ and $u(i)v$ such that 
$$\begin{array}{rl}
\be(|\be|\!-\!1\!+\!r\!-\!s)=&v(|v|\!-\!1\!-\!s)u(|u|\!-\!1\!+\!r)
+\dsp{\sum_{i\geq m}}a_i 
v(|v|\!-\!1\!-\!s\!+\!r\!-\!i)u(|u|\!-\!1\!+\!i) \cr
   &\mbox{}\quad +\dsp{\sum_{i\geq m}}b_i 
u(|u|\!-\!1\!-\!s\!+\!r\!-\!i)v(|v|\!-\!1\!+\!i)
\end{array}$$
with $a_i, b_i\in \C$. In particular, 
when $s=r=0$, then $\be=v\ast_mu$ is the product in $A_m(V)$. 
We will use this argument in several 
places. 

Set $D_n=\left< \al\in V \mid 
\mbox{ $\al=v(|v|\!-\!1\!+\!n)u(|u|\!-\!1\!-\!n)$ 
in $A(V)$ for 
some $v,u\in V$ } \right>$. \\
Clearly, $D_n$ contains $O(V)$ and $D_n/O(V)$ is an ideal 
of $A(V)$. 
If $W$ is an irreducible $V$-module and $0\not=w\in W(0)$, 
then there is an element 
$v\in V$ such that $0\not=v(|v|\!-\!1\!-\!n)w\in W(n)$ 
since $W(n)\not=0$. There is also an element $u\in V$ such that 
$u(|u|\!-\!1\!+\!n)v(|v|\!-\!1\!-\!n)w\not=0$, (see \cite{DM}).  
Hence $D_n$ covers $\End(W(0))$ for all simple modules $W$ 
and so $A(V)=D_n+J(A(V))$, which means $D_n=V$. 
\prend

The purpose in this subsection is to prove the following three 
propositions,  
which we will use in the next section.

\begin{prn}
If $S(v,\tau)=\dsp{\sum_{i=0}^{\infty}} S_i(v)q^i\in \C[[q]]$ 
satisfies that 
$S(\al,\tau)=0$ for all $\al\in O_q(V)$, then 
$S(v,\tau)\in \C[[q]]q^{n+1}$ for all $v\in O_n(V)$. 
\end{prn}

In particular, $S_i$ is a symmetric linear function of $A_i(V)$. 

\begin{prn}
Assume that $S(v,\tau)
=\dsp{\sum_{j=0}^N}\left(\dsp{\sum_{i=0}^{\infty}}
S_{ji}(v)q^{i+r}\right)(2\pi i\tau)^j
\in \C[[q]]q^r[\tau]$ 
satisfies $S(\al,\tau)=0$ for all $\al\in O_q(V)$ and 
$$S(L[\!-\!2]v-\sum_{k=1}^{\infty} E_{2k}(\tau)L[2k\!-\!2]v,\tau)
=\frac{1}{2\pi i}\frac{d}{d\tau}S(v,\tau)$$ 
for all $v\in V$. 
Then $S_{jm}
\left( (\omega\!-\!r\!-\!\frac{c}{24}\!-\!m)^{N-j+1}
\ast_mv\right)=0$ for any $m$ and $j$. 
\end{prn}

\vspace{1mm}

\begin{prn}
Assume that $S(v,\tau)
=\dsp{\sum_{i=0}^{\infty}} S_i(v)q^i\in \C[[q]]$ 
satisfies 
$S(\al,\tau)=0$ for all $\al\in O_q(V)$ and 
$S_m((\omega\!-\!r\!-\!\frac{c}{24}-m)^s\ast_mv)=0$ 
for all $v\in V$ and $m$. 
If $S_n=0$ for some $n>{\it l}\!+\!r$, then $S_0=0$, where 
${\it l}$ is given in 
the beginning of this section. 
\end{prn}

\noindent
{\bf Proof of Proposition 4.4, 4.5 and 4.6} \\
In order to prove the three propositions above at once, we will 
review the proof of Proposition 4.3.3 in \cite{Zh}. 
Zhu first obtained 
$$\begin{array}{l}
\dsp{\tr_{|M} v(|v|\!-\!1\!-\!k)u(|u|\!-\!1\!+\!k)q^{L(0)}
=\tr_{|M}
\frac{-q^k}{1-q^k}\sum_{s=0}^{\infty}\binom{|v|\!-\!1\!-\!k}{s}
o(v(s)u)q^{L(0)}}  \cr
\dsp{\tr_{|M} u(|u|\!-\!1\!-\!k)v(|v|\!-\!1\!+\!k)q^{L(0)}
=\tr_{|M}
\frac{q^k}{1-q^k}\sum_{s=0}^{\infty}\binom{|v|\!-\!1\!+\!k}{s}
o(v(s)u)q^{L(0)}},
\end{array}$$ 
for $V$-modules $M$ and $k\not=0$, 
see the proof of Proposition 4.3.2 in \cite{Zh}. 

Using the equation above, Zhu denoted 
$$\tr_{|M} w^{|v|}z^{|u|}Y(v,w)Y(u,z)q^{L(0)}\mbox{  and  }
\tr_{|M} z^{|u|}w^{|v|}Y(u,z)Y(v,w)q^{L(0)} $$
as an infinite linear combinations of $o(v)o(u)$, 
$o(v[0]u)=o(v)o(u)\!-\!o(u)o(v)$, \\
$\dsp{\tr_{|M}\frac{-q^k}{1\!-\!q^k}\sum_{s=0}^{\infty}
\binom{|v|\!-\!1\!-\!k}{s}o(v(s)u)q^{L(0)}(\frac{w}{z})^k}$ 
and \\
$\dsp{\tr_{|M}\frac{q^k}{1\!-\!q^k}
\sum_{s=0}^{\infty}\binom{|v|\!-\!1\!+\!k}{s}
o(v(s)u)q^{L(0)}(\frac{z}{w})^k}$. \\
Then substituting them into the normal product 
$$Y(v(i)u,z)=
\Res_w\{ (w\!-\!z)^iY(v,w)Y(u,z)-(\!-\!z\!+\!w)^iY(u,z)Y(v,w)\}$$
and using an expansion 
$v[\!-\!1]u=\sum_{i\geq -1}c_i v(i)u$ with $c_i\in \C$, 
he expressed the term 
$o(v[\!-\!1]u)$ between $\tr_{|M}$ and $q^{L(0)}$ 
as a linear combinations of \\
$o(v)o(u)$, $o(v[0]u)$, $\frac{-q^k}{1-q^k}
\sum_{s=0}^{\infty}\binom{|v|\!-\!1\!-\!k}{s}o(v(s)u)$ 
and $\frac{q^k}{1-q^k}
\sum_{s=0}^{\infty}\binom{|v|\!-\!1\!+\!k}{s}o(v(s)u)$.

The next step is the crucial part of his paper.  He changed the 
shape of the above expression of $o(v[\!-\!1]u)$ into 
$$\underbrace{o(v)o(u)\!-\!o(v[0]u)}_{o(u)o(v)}
+\sum_{k=1}^{\infty}E_{2k}(q)o(v[2k\!-\!1]u)   $$
using the equations (c.f. (4.3.8)-(4.3.11) and Proposition 4.3.2 
in \cite{Zh}). Namely, he proved 
$$\tr_{|M}\left\{ o(v[\!-\!1]u)\!-\!o(v)o(u)\!+\!o(v[0]u)
-\sum_{k=1}^{\infty}E_{2k}(q)o(v[2k\!-\!1]u)\right\}q^{L(0)}=0. $$
He obtained two important equations from it. 
The first is, by substituting $\tilde{\omega}$ for $v$, 
$$\tr_{|M}\left\{o(L[\!-\!2]u)
-\sum_{k=1}^{\infty}E_{2k}(\tau)o(L[2k\!-\!2]u)\right\}
q^{L(0)-c/24}=\frac{1}{2\pi i}\frac{d}{d\tau} 
\tr_{|M}o(u)q^{L(0)-c/24}. $$
The second is, by substituting $L[\!-\!1]v$ for $v$, he had 
$$ \tr_{|M}\left\{ 
v[\!-\!2]u+\sum_{k=2}^{\infty}(2k\!-\!1)v[2k\!-\!2]u
E_{2k}(q)\right\}q^{L(0)}=0.$$
These equations are so beautiful that we can see the modular 
invariance property of $O_q(V)$ clearly. Once we know the 
modular invariance 
property of $O_q(V)$, we don't need these forms. We are 
interested in 
the structure of $O_q(V)$ from a view point of ordinary 
vertex operators $Y(v,z)=\sum_{i\in \Z} v(i)z^{-i-1}$. 
We go back to the former form of 
$$o(u)o(v)+\sum_{k=1}^{\infty}E_{2k}(\tau)o(v[2k\!-\!1]u). $$
Namely, we change an expansion of $o(v[\!-\!1]u)$: 
$$\begin{array}{l}
\dsp{o(v[\!-\!1]u)=\sum_{i\geq -1} c_i o(v(i)u) }\cr
\mbox{}\quad\dsp{=\sum_{i\geq -1}c_i
\sum_{j=0}^{\infty}(\!-\!1)^j\binom{i}{j}
\left\{
v(i\!-\!j)u(|v|\!+\!|u|\!-\!i\!-\!2\!-\!j)
-(\!-\!1)^iu(|v|\!+\!|u|\!-\!2\!-\!j)v(j)\right\} }
\end{array}$$
by replacing 
$v(|v|\!-\!1\!-\!k)u(|u|\!-\!1\!+\!k)$ and 
$u(|u|\!-\!1\!-\!k)v(|v|\!-\!1\!+\!k)$ by 
$\dsp{\frac{-q^k}{1\!-\!q^k}
\sum_{s=0}^{\infty}\binom{|v|\!-\!1\!-\!k}{s}v(s)u}$ and 
$\dsp{\frac{q^k}{1\!-\!q^k}
\sum_{s=0}^{\infty}\binom{|v|\!-\!1\!+\!k}{s}v(s)u}$, 
respectively, for $k\not=0$.   
We note $c_{-1}=1$. 

To simplify the arguments, we express this process by 
notation $\theta$, that is, for $\al=\sum_i \la_i v(i)u$, 
we first develop $ \sum_i \la_i o(v(i)u)$ as a linear sum 
$$a o(v)o(u)+bo(u)o(v)+\sum_{k\not=0} \la_k 
v(|v|\!-\!1\!-\!k)u(|u|\!-\!1\!+\!k)
+\sum_{k\not=0} \mu_k u(|u|\!-\!1\!-\!k)v(|v|\!-\!1\!+\!k)$$
with $a, b, \la_i, \mu_i\in \C$ 
by using associativity, then define 
$$ \begin{array}{l}
\dsp{\theta(\sum_i \la_i v(i)u\!-\!ao(v)o(u)\!-\!bo(u)o(v))}\cr
\mbox{}\quad \dsp{=\sum_{k\not=0} \la_k \frac{-q^k}{1-q^k}
\sum_{s=0}^{\infty}\binom{|v|\!-\!1\!-\!k}{s}v(s)u
+\sum_{k\not=0}\mu_k \frac{q^k}{1-q^k}
\sum_{s=0}^{\infty}\binom{|v|\!-\!1\!+\!k}{s}v(s)u}.
\end{array}$$
What Zhu has obtained are 
$$
v[\!-\!1]u-\theta(v[\!-\!1]u\!-\!o(u)o(v))=v[\!-\!1]u
-\sum_{k=1}^{\infty}E_{2k}(\tau)v[2k\!-\!1]u 
\eqno{(4.3)} $$
and 
$$v[\!-\!2]u-\theta(v[\!-\!2]u)=v[\!-\!2]u
+\sum_{k=1}^{\infty} (2k\!-\!1)E_{2k}(\tau)v[2k\!-\!2]u
\in O_q(V) \eqno{(4.4)}$$ 

On the other hand, if $i\geq 0$, then we have 
$$ \begin{array}{l}
o(v(i)u)\cr
\mbox{}\qquad \dsp{=\sum_{j=0}^{i}(\!-\!1)^j
\binom{i}{j}\{v(i\!-\!j)u(|v|\!+\!|u|\!-\!i\!-\!2\!+\!j)
-(\!-\!1)^iu(|v|\!+\!|u|\!-\!2\!-\!j)v(j)\}}\cr
\mbox{}\qquad \dsp{=\sum_{j=0}^{i}(\!-\!1)^{i+j}\binom{i}{j}
\{v(j)u(|v|\!+\!|u|\!-\!2\!+\!j)
-u(|v|\!+\!|u|\!-\!2\!-\!j)v(j)\}}.
\end{array}$$
If we replace  
$ v(|v|\!-\!1\!+\!k)u(|u|\!-\!1\!-\!k)
-u(|u|\!-\!1\!-\!k)v(|v|\!-\!1\!+\!k)$ by  
$$\frac{1}{1\!-\!q^k}
\sum_{s=0}^{\infty}\binom{|v|\!-\!1\!+\!k}{s}v(s)u
-\frac{q^k}{1\!-\!q^k}
\sum_{s=0}^{\infty}\binom{|v|\!-\!1\!+\!k}{s}v(s)u 
=\sum_{s=0}^{\infty}\binom{|v|\!-\!1\!-\!k}{s}v(s)u$$
for each $k\not=0$, then 
$$o(v(i)u)=\sum_{j=0}^{i}(\!-\!1)^{i+j}\binom{i}{j}
\{v(j)u(|v|\!+\!|u|\!-\!2\!+\!j)
-u(|v|\!+\!|u|\!-\!2\!-\!j)v(j)\}$$
is replaced by 
$$\sum_{j=0}^i(\!-\!1)^{i+j}\binom{i}{j}
\sum_{s=0}^{\infty}\binom{j}{s}v(s)u=v(i)u . $$
Namely, $\theta(v(i)u)=v(i)u$ for $i\geq 0$.
By cancelling $v(i)u$  $(i\geq 0)$ from both sides of $v[\!-\!1]u$ 
and $\theta(v[\!-\!1]u$ in (4.3) and also both sides of 
$v[\!-\!2]u$ and $\theta(v[\!-\!2]u)$ in (4.4),   
we obtain 
$$
v(\!-\!1)u-\theta\left(v(\!-\!1)u-o(u)o(v)\right)
=v[\!-\!1]u-\sum_{k=1}^{\infty}E_{2k}(\tau)v[2k\!-\!1]u 
\eqno{(4.5)} $$
$$  \omega(\!-\!1)u-\frac{c}{24}u-
\theta(\omega(\!-\!1)u\!-\!o(u)o(\omega))
=L[\!-\!2]u-\sum_{k=1}^{\infty}E_{2k}(\tau)L[2k\!-\!2]u, 
\quad \mbox{ and} \eqno{(4.6)} $$
$$\begin{array}{l}
\dsp{v(\!-\!2)u+|v|v(\!-\!1)u-\theta(v(\!-\!2)u)+|v|v(\!-\!1)u)}\cr
\mbox{}\qquad \qquad\dsp{=v[\!-\!2]u
+\sum_{k=1}^{\infty}
(2k\!-\!1)E_{2k}(\tau)v[2k\!-\!2]u\in O_q(V)}, 
\end{array}\eqno{(4.7)}$$ 
since $\theta({\bf 1}(\!-\!1)u-o(u))=0$. 
Substituting $L(\!-\!1)^mv$ into $v$ of (4.7), we have 
$$v(\!-\!2\!-\!m)u+\frac{|v|\!+\!m}{m+1}v(\!-\!1\!-\!m)u-
\theta(v(\!-\!2\!-\!m)u+\frac{|v|\!+\!m}{m+1}v(\!-\!1\!-\!m)u) 
\in O_q(V).$$

If $a\in O_n(V)$, then $a$ is a linear combination of 
$$v\circ_nu=\Res_z Y(v,z)u\frac{(1+z)^{|v|+n}}{z^{2+2n}}$$
for some $v,u\in V$ and the expansion of 
$v\circ_nu$ by associativity is an 
infinite linear sums of 
$v(|v|\!-\!1\!-\!k)u(|u|\!-\!1\!+\!k)$ and 
$u(|u|\!-\!1\!-\!k)v(|v|\!-\!1\!+\!k)$ with 
$k\geq n\!+\!1$. 
Since for $k\geq n\!+\!1$, 
$$ \theta\left(v(|v|\!-\!1\!-\!k)u(|u|\!-\!1\!+\!k)\right)
=\frac{-q^k}{1-q^k}
\sum_{s=0}^{\infty}\binom{|v|\!-\!1\!-\!k}{s}v(s)u $$
and 
$$\theta\left(u(|u|\!-\!1\!-\!k)v(|v|\!-\!1\!+\!k)\right)
=\frac{q^k}{1-q^k}
\sum_{s=0}^{\infty}\binom{|v|\!-\!1\!+\!k}{s}v(s)u $$
are contained in $V[q]q^{n+1}$, we obtain 
$S(a,q)\in V[[q]]q^{n+1}$, which proves Proposition 4.4. \\

By the same arguments, if 
$$S(v,\tau)=\sum_{j=0}^NS_j(v,\tau)q^r(2\pi i\tau)^j
=\sum_{j=0}^N\left(\sum_{i=0}^{\infty}S_{ji}(v)q^{i}\right)q^r
(2\pi i\tau)^j$$
satisfies $S(\al,\tau)=0$ for all $\al\in O_q(V)$ and 
$$ S(L[\!-\!2]v-\sum E_{2k}(\tau)L[2k\!-\!2]v,\tau)
=\frac{1}{2\pi i}\frac{d}{d\tau}S(v,\tau) $$
for all $v\in V$, then 
$$S(\sum_{i>-M} \la_i \omega(i)u,\tau)=s\left(
\frac{1}{2\pi i}\frac{d}{d\tau}S(u,\tau)
+\frac{c}{24}S(u,\tau)\right)
+S(\theta(\sum_{i>-M} \la_i\omega(i)u-so(\omega)o(u)),\tau) $$ 
for $\sum_{i>-M} \la_i\omega(i)u$, where $s$ is given by 
$$o(\sum_{i>-M} \la_i\omega(i)u)=s~o(\omega)o(u)
+\sum_{j>0}\left(a_j\omega(1\!-\!j)u(|u|\!-\!1\!+\!j)
+b_ju(|u|\!-\!1\!-\!j)\omega(1\!+\!j)\right)$$
with $s, a_j, b_j\in \C$.  

It follows from the definition of $\ast_n$ that 
$o(\omega\ast_n u)\in \omega(1)o(u)+O_n(V)$ and so 
$$S(\omega\ast_nu,\tau)=\frac{1}{2\pi i}\frac{d}{d\tau}S(u,\tau)
+\frac{c}{24}S(u,\tau)
+S(\theta(\omega\ast_nu\!-\!o(\omega)o(u)),\tau). $$
Therefore we obtain 
$$\begin{array}{l}
\dsp{\sum_{j=0}^N\sum_{i=0}^{\infty}S_{ji}(\omega\ast_nu)q^{i+r}
(2\pi i\tau)^j
=\sum_{j=0}^NS_j(\omega\ast_nu,\tau)q^r(2\pi i\tau)^j
=S(\omega\ast_nu,\tau) }\cr
\dsp{=\frac{1}{2\pi i}\frac{d}{d\tau}S(u,\tau)
+\frac{c}{24}S(u,\tau)
+S(\theta(\omega\ast_nu\!-\!o(\omega)o(u)),\tau) }\cr
\dsp{=\frac{1}{2\pi i}\frac{d}{d\tau}\left(
\sum_{j=0}^NS_{ij}(u)q^{i+r}(2\pi i\tau)^j\right)
+\frac{c}{24}\sum_{j=0}^NS_{ij}(u)q^{i+r}(2\pi i\tau)^j }\cr
\mbox{}\quad 
\dsp{+\sum_{j=0}^NS_j(\theta(\omega\ast_nu\!-\!o(\omega)o(u)),\tau)
(2\pi i\tau)^j}\cr
\dsp{=\sum_{j=0}^N\sum_{i=0}^{\infty}
S_{ij}(u)(i+r)q^{i+r}(2\pi i\tau)^j
+\sum_{j=0}^N\sum_{i=0}^{\infty}
S_{ij}(u)q^{i+r}(2\pi i\tau)^{j-1}}\cr
\mbox{}\quad\dsp{ +\frac{c}{24}\sum_{j=0}^N
\sum_{i=0}^{\infty}S_{ij}(u)q^{i+r}
(2\pi i\tau)^j 
+\sum_{j=0}^NS_j(\theta(\omega\ast_nu\!-\!o(\omega)o(u)),\tau)
(2\pi i\tau)^j}\cr
\dsp{=\sum_{j=0}^N\sum_{i=0}^{\infty}
S_{ij}((i+r)u)q^{i+r}(2\pi i\tau)^j
+\sum_{j=0}^N\sum_{i=0}^{\infty}
S_{ij}(u)q^{i+r}(2\pi i\tau)^{j-1}}\cr
\mbox{}\quad \dsp{+\sum_{j=0}^NS_{ij}(\frac{c}{24}u)q^{i+r}
(2\pi i\tau)^j
+\sum_{j=0}^NS_j(\theta(\omega\ast_nu\!-\!o(\omega)o(u)),\tau)
(2\pi i\tau)^j}
\end{array}$$
and so 
$$\begin{array}{l}
\dsp{\sum_{i=0}^{\infty}
S_{j,i}(\omega\ast_nu)q^{i+r}(2\pi i\tau)^j }\cr
\dsp{=\sum_{i=0}^{\infty}S_{j,i}((i\!+\!r\!-\!\frac{c}{24})u)
q^{i+r}(2\pi i\tau)^j
+\sum_{i=0}^{\infty}S_{j+1,i}(u)q^{i+r}(2\pi i\tau)^j }\cr
\mbox{}\quad 
+\dsp{S_j(\theta(\omega\ast_nu\!-\!o(\omega)o(u)),\tau)
(2\pi i\tau)^j}
\end{array}$$
for each $j$. 
Thus we have 
$$S_{j,i}((\omega\!-\!i\!-\!r\!-\!\frac{c}{24})\ast_nu)
=S_{j+1,i}(u) $$
for $i\leq n$ since 
$\theta(\omega\ast_nu\!-\!o(\omega)o(u))\in V[[q]]q^{n+1}$.
It follows from $S_{N+1,n}(u)=0$ that \\
$S_{j,n}((\omega\!-\!n\!-\!r\!-\!\frac{c}{24})^{N-j+1}\ast_n u)
=0$, which proves Proposition 4.5.\\

Suppose that Proposition 4.6 is false. Namely, 
there is an integer $n>{\it l}+r$ such that $S_n=0$ 
and $S_0(\al)\not=0$ for some $\al\in V$. 
By Proposition 4.4, $S_0$ is a (symmetric) linear map of 
$A(V)=V/O(V)$. Set $A=A(V)/\Rad(S_0)$. Then 
$(\omega\!-\!r\!-\!\frac{c}{24})^sA=0$. 
By Lemma 4.3, there are $v^i,u^i\in V$ such that 
$\sum_{i=1}^p v^i(|v^i|\!-\!1\!+\!n)u^i(|u^i|\!-\!1\!-\!n)
=o(\al)$ in $A(V)$. 
By the choice of $\al$, we may assume 
$v(|v|\!-\!1\!+\!n)u(|u|\!-\!1\!-\!n)=o(\al)$ in $A(V)$ 
for some $v,u\in V$.
 
As we mentioned in the proof of Lemma 4.3, 
there is an element $\be\in V$ such that 
$$\begin{array}{rl}
o(\be)=&u(|u|\!-\!1\!-\!n)v(|v|\!-\!1\!+\!n)
+\dsp{\sum_{i>n}} a_iv(|v|\!-\!1\!-\!i)u(|u|\!-\!1\!+\!i)\cr
 &\mbox{}\quad +\dsp{\sum_{i>n}} 
b_iu(|u|\!-\!1\!-\!i)v(|v|\!-\!1\!+\!i) 
\end{array}$$ 
for some $a_i, b_i\in \C$. Then we obtain 
$$\begin{array}{rl}
o(\al)=&\dsp{v(|v|\!-\!1\!+\!n)u(|u|\!-\!1\!-\!n)
=[v(|v|\!-\!1\!+\!n),u(|u|\!-\!1\!-\!n)] }\cr
=&\dsp{\sum_{i=0}^{\infty}\binom{|v|\!-\!1\!+\!n}{i}o(v(i)u)}
\end{array}$$ 
on $W(0)$ for any $\N$-graded weak $V$-modules $W$ and so 
$\al=\dsp{\sum_{i=0}^{\infty}}\binom{|v|\!-\!1\!+\!n}{i}v(i)u$ 
in $A(V)$. On the other hand, since $n>1$, 
we have $\be\in O(V)$ and 
$$\begin{array}{rl}
S(\be,\tau)=&S(\theta(\be),\tau)\cr
\in &\dsp{\frac{q^n}{1-q^n}S(\sum_{i=0}^{\infty}
\binom{|v|\!-\!1\!+\!n}{i}v(i)u,\tau)q^n+q^{n+1}\C[[q]]}\cr
&\dsp{=S(\sum_{i=0}^{\infty}\binom{|v|\!-\!1\!+\!n}{i}v(i)u,\tau)q^n
+q^{n+1}\C[[q]]=S(\al,\tau)q^n+q^{n+1}\C[[q]]}. 
\end{array}$$
Since coefficients of $S(\be,\tau)$ at $q^n$ are always zero 
and the constant 
term of $S(\al,\tau)$ is nonzero, we have a contradiction. 

This completes the proofs of three propositions. 
\prend

\section{The space of one point functions on the torus}
In this section, we will just follow the proofs in \cite{Zh} and 
\cite{DLiM3} with suitable modification since 
we use pseudo-trace functions which satisfy the same properties 
as do the ordinary trace functions and so 
we will skip the most part 
of the proof. The differences between our case and Zhu's case 
(and also the case in \cite{DLiM3}) are that we will treat 
$A_n(V)$ and our $A_n(V)$ is not a semisimple algebra, 
for example, $\omega$ might not act on $A_n(V)$ semisimply. 
However, since $A_n(V)$ is a finite dimensional algebra, 
there are $r_i\in \C$ and $\mu(r_i)\in \Z$ such that 
$\prod_{i=1}^s
(\omega\!-\!\frac{c}{24}\!-\!r_i\!-\!n)^{\mu(r_i)}A_n(V)=0$. \\

Let's recall the following notation from \cite{Zh} 
and \cite{DLiM3}: Consider $V[E_4(q),E_6(q)]\subseteq V[[q]]$. 
$O_q(V)$ is the submodule of $V[E_4(q),E_6(q)]$ generated by 
elements of the type 
$$  v[\!-\!2]u+ \sum_{k=2}^{\infty} (2k\!-\!1)E_{2k}(\tau)
\otimes v[2k\!-\!2]u  \quad \mbox{ with } \quad v,u \in V $$
and 
$$ v[0]u  \quad \mbox{ with } \quad v,u\in V. $$

We also recall the definition of the space $\CC_1(V)=\CC(1,1)$ of 
one point functions with trivial automorphisms from \cite{DLiM3}. 

\begin{dfn}
We define the space ${\CC}_1(V)$ of one point functions on $V$ 
to be the $ \C$-linear space consisting of functions
$$  S: V[E_4(q),E_6(q)]\otimes \CH \to \C   $$
satisfying the following conditions: \\
{\rm (C1)}  For $u \in V(\Gamma(1))$, $S(u,\tau)$ is holomorphic 
in $ \tau\in \CH$. \\
{\rm (C2)}  $S(\sum f_i(\tau)\otimes u_i, \tau)
=\sum_i f_i(\tau)S(u_i, \tau)$ 
for $f_i(\tau) \in \C[E_4(q),E_6(q)]$ and $u_i \in V$. \\
{\rm (C3)}  For $u \in O_q(V)$, $S(u, \tau)=0$.  \\
{\rm (C4)}  For $u \in V$, 
$$  S(L[\!-\!2]u, \tau)=\frac{1}{2\pi i}\frac{d}{d\tau}S(u, \tau)+
\sum_{k=1}^{\infty}E_{2k}(\tau)S(L[2k\!-\!2]u, \tau).$$
\end{dfn}

By exactly the same proof, we have the following modular 
invariance property of $\CC_1(V)$, (see Theorem 5.1.1 in \cite{Zh} 
and Theorem 5.4 in \cite{DLiM3}.)

\begin{thm}(Modular-Invariance). For $S\in \CC_1(V)$ and 
$\ga\in \binom{a\quad b}{c\quad d}\in SL(2,\Z)$, define 
$S|\ga(v,\tau)=\frac{1}{(c\tau+d)^h}S(v,\ga(\tau))$ 
for $v\in V_{[h]}$ and extend linearly. Then $S|\ga\in \CC_1(V)$. 
\end{thm}

Set 
$$S^W(u, \tau)=\tr^{\phi}_{W}o(u)q^{L(0)-c/24}$$ 
for a generalized Verma module $W$ interlocked with 
a symmetric linear function $\phi$ of $A_n(V)$ 
and an element $u\in V$, 
where $c$ is the central charge of $V$. 
Then extend it linearly for $V[E_4(q),E_6(q)]$,

We next prove that $S^W(u,\tau)\in \CC_1(V)$. What we have to do 
is to show that $S^W(u,\tau)$ is holomorphic on the 
upper half plane. 
In order to prove this fact, we will show that $S^W$ and 
all $S\in \CC_1(V)$ satisfy differential equations.
The proof is essentially the same as the arguments (4.4.11) in 
\cite{Zh} and Lemma 6.1 in \cite{DLiM3}. 

\begin{thm} Assume that $V$ is $C_2$-cofinite. 
If $S(v,\tau)=0$ for $v\in O_q(V)$ and 
$$ S(L[\!-\!2]u-\sum_{k=1}^{\infty}E_{2k}(\tau)L[2k\!-\!2]u,\tau)
=\frac{1}{2\pi i\tau}\frac{d}{d\tau}S(u,\tau)$$
for $u\in V$, then 
there are $m\in \N$ and 
$r_i(\tau)\in \C[E_4(q),E_6(q)]$, $0\leq i\leq m\!-\!1$, 
such that 
$$S(L[\!-\!2]^mv,\tau)+
\sum_{i=0}^{m-1}r_i(\tau)S(L[\!-\!2]^iv,\tau)=0. $$
In particular, $S(u, \tau)$ converges absolutely and uniformly 
in every closed subset of the domain $ \{q\  \mid \ |q|<1\}$ 
for every $u \in V$ and the limit function can be written as a 
linear sum of $q^hf(q)$, where $f(q)$ is some analytic function in 
$ \{q \mid |q|<1\}$. In particular, $S\in \CC_1(V)$. 
\end{thm}

As a corollary, we obtain:

\begin{cry}
$S^W(u,\tau)$ is holomorphic on the upper half-plane for $u\in V$ 
and \\
$S^W(\ast,\tau)\in \CC_1(V)$. 
\end{cry}

\subsection{Spanning set of $\CC_1(V)$}
The remaining is to show that $\CC_1(V)$ is spanned 
by pseudo-trace functions $S^W(\ast,\tau)$. 
So we will proved the following main theorem, which 
covers a nonsemisimple version 
of Theorem 5.3.1 in \cite{Zh}. 

\begin{thm}  Suppose that $V=\oplus_{m=0}^{\infty}V_m$ 
is a $C_{2}$-cofinite VOA with 
central charge $c$. Take an integer $n$ 
sufficiently large and let 
$ \{W^1, ..., W^m\}$ be the set of $n$-th generalized 
Verma $V$-modules $W^i$ 
interlocked with some symmetric linear function $\phi^i$ 
of $A_n(V)$.  
Then ${\CC}_1(V)$ is spanned by 
$$ \{ S^{W^1}(\cdot, \tau),..., S^{W^m}(\cdot,\tau) \}.  $$
In particular, $\dim \CC_1(V)\!=\!
\dim A_n(V)/[A_n(V),A_n(V)]\!
-\!\dim A_{n-1}(V)/[A_{n-1}(V),A_{n-1}(V)]$.
\end{thm}

We note that the dimension of $A_n(V)$ is finite and 
so $\CC_1(V)$ is also of 
finite dimension. \\

\noindent
{\bf Proof of Theorem 5.5} 

Let $S\in \CC_1(V)$. 
We will prove that $S$ is a sum of pseudo-trace 
functions. 
By the same arguments as in \cite{Zh}, it follows 
from Theorem 5.3 that 
there are integers $d, N_1,...,N_d$ which 
does not depend on $v$ such that 
$$ S(v,\tau)=\sum_{s=0}^d S_s(v,\tau)q^{r_s}$$
and each $S_s(v,\tau)$ can be further decomposed as 
$$S_s(v,\tau)=\sum_{j=0}^{N_s}S_{sj}(v,\tau)(2\pi i\tau)^j,$$ 
where $r_{1},...,r_{d}$ are complex 
numbers independent of $v$, $r_{s_1}\!-\!r_{s_2}\not\in \Z$ for 
$s_1\not=s_2$ and $S_{sj}(v,\tau)$ 
has a $q$-expansion $S_{sj}(v,\tau)
=\sum_{i=0}^{\infty}C_{sji}(v)q^k$ 
with $C_{sji}(v)\in \C$ and 
for each $j$ there is $s$ such that $C_{sj0}\not=0$. 
Since $r_{s_1}\!-\!r_{s_2}\not\in \Z$ for 
$s_1\not=s_2$, each 
$S_s(v,\tau)q^{r_s}$ satisfies (C2)$\sim$(C4) and also (C1) by 
Theorem 5.3. Hence we may assume $d=1$ and  
$$ S(v,\tau)=\sum_{j=0}^N S_j(v,\tau)(2\pi i\tau)^j=
\sum_{j=0}^N
\left(\sum_{i=0}^{\infty}C_{ji}(v)q^{i+r}\right)(2\pi i\tau)^j, 
\eqno{(5.1)}$$
with $r\in \C$.  
In the case where $A(V)$ is semisimple, they proved 
in \cite{Zh} and \cite{DLiM3} 
that $N=0$ and 
$\sum_{i=0}^{\infty}S_{i}(v)q^{i+r}$ is a sum of trace functions. 
However, if $A_n(V)$ is not semisimple, 
we may have nonzero $N$ since we consider logarithmic 
modules, too. In order to continue the proof, 
we need the following two lemmas.

\begin{lmm} 
Let $W=W^{(n)}_T$ be an $n$-th generalized Verma module 
interlocked with $(P,\phi)$ satisfying 
$(\omega\!-\!r\!-\!n\!-\!\frac{c}{24})^{N+1}W^{(n)}_T(n)=0$ 
for some $r\in \C$ and $N\in \N$. 
Then there are 
constants $b_1,...,b_{s-1}$ such that 
$$\tr^{\phi}_W o(v)q^{L^s(0)-c/24}=S^W(v,\tau)
-\sum_{i=1}^{N}b_iS^{W/W\FN^i}(v,\tau)(2\pi i\tau)^i,$$ 
where $\FN^i
=\{a\in P \mid (\omega\!-\!r\!-\!n\!-\!\frac{c}{24})^ia=0\}$ and 
$L^s(0)$ is a semisimple part of $L(0)$, which acts on $W(m)$ as 
$m\!+\!r\!+\!\frac{c}{24}$ and $(L(0)\!-\!L^s(0))^{N+1}=0$ on $W$. 
\end{lmm}

\pr 
We first note that there are constants $b_0=1,b_1,...$ such that 
$e^{2\pi i\al}
=1\!+\! b_1\al e^{2\pi i\al}\!+\!b_2\al^2e^{2\pi i\al}
\!+\!\cdots\!+\!b_{s-1}\al^{s-1}e^{2\pi i\al}\!+\!\cdots$. 
Hence we obtain 
$$\begin{array}{rl}
S^W(v,\tau)=&\dsp{\tr^{\phi} o(v)q^{L(0)-\frac{c}{24}}
=\tr^{\phi} o(v)q^{(L^s(0)-\frac{c}{24})}q^{(L(0)-L^s(0))} }
\cr
=&\dsp{\tr^{\phi}_{W} o(v)q^{L^s(0)-\frac{c}{24}}
\!+\!\sum_{i=1}^{N}\tr^{\phi}_Wo(v)b_i
(L(0)\!-\!L^s(0))^i(2\pi i\tau)^i 
q^{L^s(0)-c/24+(L(0)-L^s(0))} }\cr
=&\!\!\dsp{\tr^{\phi}_{W} o(v)q^{L^s(0)-\frac{c}{24}}
\!+\!\sum_{i=1}^{N}\tr^{\phi}_Wo(v)
b_i(\omega\!-\!n\!-\!\frac{c}{24}\!-\!r)_P^i
q^{L(0)-\frac{c}{24}}(2\pi i\tau)^i }\cr
=&\!\!\!\dsp{\tr^{\phi}_W o(v)q^{L^s(0)-\frac{c}{24}}\!
+\!\sum_{i=1}^{N}
b_i\left(\tr^{(\omega-n-\frac{c}{24}-r)^i\phi}_{W/W\FN^i}o(v)
q^{L(0)-\frac{c}{24}}(2\pi i\tau)^i\right)} \qquad 
\mbox{ by (3.6)}\cr
=&\!\!\!\dsp{S^W(v,\tau)+\sum_{i=1}^{N}
b_iS^{W/W\FN^i}(v,\tau)(2\pi i\tau)^i}. \prend
\end{array}$$

We will prove that $S_{0}(v,\tau)
=\sum_{i=0}^{\infty} C_{0i}(v)q^{i+r}$ is a linear sum 
of pseudo-trace functions 
with semisimple grading operator $L^s(0)$, say 
$S_0(v,\tau)=\sum_p a_p \tr^{\phi^p}_{W^p} o(v)q^{L^s(0)-c/24}.$ 
Then 
$$\tilde{S}(v,\tau)=S(v,\tau)
-\sum_p a_p\left(S^{W^p}(v,\tau)
-\sum_{i=1}^{s-1}b_iS^{W/W\FN^i}(v,\tau)
(2\pi i\tau)^i\right) \in \CC_1(V),$$ 
but if we express it by 
$$\tilde{S}(v,\tau)=\sum_{s=0}^N \sum_{j=0}^{\infty}
\tilde{S}_{sj}(v)q^{j+r}(2\pi i\tau)^s,$$ 
then $\tilde{S}_{0j}(v)=0$ for all $j$ and $v\in V$. 
However, since 
$\tilde{S}(L[\!-\!2]u
-\sum_{k=1}^{\infty}E_{2k}(\tau)L[2k\!-\!2]u,\tau)
=\frac{1}{2\pi i}\frac{d}{d\tau}\tilde{S}(u,\tau)$, 
$S_{0,j}=0$ implies $S_{1,j}=0$ and so on. 
We hence have $\tilde{S}(v,\tau)=0$ for all $v\in V$ as desired. 

So it is sufficient to prove that 
$S_{0}(v,\tau)=\sum_{i=0}^{\infty} C_{0i}(v)q^{i+r}$ is 
the coefficient of $(2\pi i\tau)^0$-term of a linear sum 
of pseudo-trace functions $S^{W^p}(v,\tau)$. 

Since 
$S(L[\!-\!2]u-\sum_{k=1}^{\infty}E_{2k}(\tau)L[2k\!-\!2]u,\tau)
=\frac{1}{2\pi i}\frac{d}{d\tau}S(v,\tau)$, we obtain the 
following lemma from 
Proposition 4.4 and 4.5. 

\begin{lmm} We have 
$$\begin{array}{l}
C_{0,n}(v\ast_n u)=C_{0,n}(u\ast_n v) \cr
C_{0,n}(v)=0  \mbox{  for  }v\in O_{n+1}(V) \cr
C_{0,n}((\omega\!-\!\frac{c}{24}\!-\!r_{sj}\!-\!n)^{N+1}\ast_n v)=0
\end{array}$$
\end{lmm}

In particular, $A_n(V)$ has a symmetric linear map $C_{0,n}$. 
Let $A_1\oplus \cdots\oplus A_k$ be the decomposition of 
$A_n(V)/J(A_n(V))$ into the direct sum of simple algebras $A_i$ and 
$\{e_i \mid i=1,...,k\}$ a set of mutually orthogonal 
primitive idempotents of $A_n(V)$ such that 
$e_i\!+\!J(A_n(V))\in A_i$. 
Set $e=\sum e_i$. 
By Theorem 3.10, there are symmetric linear functions 
$\phi_{p}$ of $A_n(V)$ 
and $A_n(V)\times eA_n(V)e$-modules $(A_n(V)/\FN_{p})e$ 
such that $(\omega\!-\!\frac{c}{24}\!-\!r\!-\!n)^{N+1}\in \FN_{p}$ 
and $C_{0,n}$ is a sum of pseudo-trace maps, that is,  
$$C_{0,n}=
\sum_p a_{p}\tr^{\phi_{p}}_{(A_n(V)/\FN_{p})e}, $$ 
where $\FN_{p}=\Rad(\phi_{p})$ and $(A_n(V)/I_p)e$ 
are all indecomposable $A_n(V)$-modules. 

We have assumed that $n$ is large enough so that 
there are no conformal 
weights greater than $r\!+\!n\!-\!{\it l}$. 
Construct $n$-th generalized Verma modules 
$$W^{p}=W^{(n)}_{T_{p}}$$ 
from $A_n(V)$-module $T_{p}=(A_n(V)/\FN_{p})\bar{e}$. 
As we showed, 
$W^{p}$ is interlocked with $\phi_{p}$ and 
$L(0)\!-\!\frac{c}{24}\!-\!r\!-\!m$ 
acts on $W^{p}(m)$ as a nilpotent operator.  
Define a pseudo-trace function 
$$ S^{W^{p}}(v,\tau)=\tr^{\phi_p} o(v)q^{L(0)-c/24}.$$
Then $\tilde{S}(v,\tau)=S(v,\tau)-\sum_p a_{p}S^{W^{p}}(v,\tau)
=\sum_{s=0}^N
(\sum_{i=0}^{\infty}\tilde{C}_{si}(v)q^{i+r})(2\pi i\tau)^s$ 
satisfies the same properties, but 
$$\tilde{C}_{0n}(v,\tau)=0 \mbox{  for all $v\in V$}. $$
Then by Proposition 4.5, we have 
$\tilde{C}_{00}=0$. Since 
$$\tilde{S}(L[\!-\!2]u\!-\!\sum_{k=1}^{\infty}
E_{2k}(\tau)L[2k\!-\!2]u,\tau)
=\frac{1}{2\pi i}\frac{d}{d\tau}\tilde{S}(u,\tau),$$ 
$\tilde{C}_{s0}(u)=0$ for all $u\in V$ 
and $s=0,1,...,N$. 
Namely, the lowest weight $\tilde{r}$ of 
$\tilde{S}(\ast,\tau)$ 
is greater than that of $S(\ast,\tau)$. Repeating these steps, 
we obtain the desired result, 
since $V$ has only finitely many lowest weights of 
pseudo-trace functions. 
This completes the proof of Theorem 5.5. 
\prend

We next consider the case $v=\1$. 
Then by Theorem 5.5, 
$$\left< \tr^{\phi}_Wq^{L(0)-c/24} \mid W
\mbox{ is interlocked with $\phi$ of $A_n(V)$}\right>$$
is $SL_2(\Z)$-invariant. 
$\tr^{\phi}_Wq^{L(0)-c/24}$ plays a role of a generalized 
character introduced in \cite{F} and so we will call it 
{\bf a generalized character} of $W$. Let's study generalized 
characters for a while. If $\phi(1)=0$ 
and $(L(0)\!-\!L^s(0))W=0$, then $\tr^{\phi}_W 1q^{L(0)}=0$ 
and so we have:

\begin{thm}
Let $V=\oplus_{m=0}^{\infty}V_m$ be a $C_2$-cofinite VOA. 
Then the space spanned by generalized characters is 
$SL_2(\Z)$-invariant. In particular, if there is no 
logarithmic module, then the space spanned by 
the set of all (ordinary) characters is $SL_2(\Z)$-invariant.  
\end{thm}

We may assume $\phi(1)=0$ and $(L(0)\!-\!L^s(0))^mW=0$. 
Let $r$ be a conformal weight of $W$.  
Then 
$$
\tr^{\phi}_W q^{L(0)-c/24}=
\tr^{\phi}_W\sum_{j=0}^m\frac{1}{j!}(L(0)\!-\!L^s(0))^j
q^{L^s(0)-c/24}(2\pi i\tau)^j\in \C[[q]]q^{r-c/24}[\tau]\tau. 
$$ 
By Lemma 4.1, $\tr^{\phi}_W(L(0)\!-\!L^s(0))^jq^{L^s(0)-c/24}
=\tr^{(\omega\!-\!n\!-\!r)^j\phi}_{W/W\FN_j}q^{L^s(0)-c/24}$ 
is a linear combination of characters. 
Therefore, we obtain:

\begin{prn} 
A generalized character is a linear combination of 
characters with coefficients in $\C[\tau]$. 
\end{prn}

As an application of Theorem 5.5, 
$\left<{\rm ch}_W(\tau)^{SL_2(\Z)}\right> $ is of finite 
dimension for an irreducible $V$-module $W$ and   
${\rm ch}_W(\frac{-1}{\tau})\in 
\sum_{i=1}^k\C[[q]][\tau]q^{r_i-c/24}$. Therefore, 
we can apply the same arguments as in 
the proofs of Proposition 3 in \cite{AM} and 
Theorem 11.3 in \cite{DLiM3} with suitable 
modifications ($q$-powers should be replaced by   
elements in $\C[[q]][\tau]$) and so 
we obtain the following corollary.

\begin{cry}
If $V=\oplus_{m=0}^{\infty}V_n$ is a $C_2$-cofinite VOA, 
then the central charge and the conformal weights 
are all rational numbers. 
\end{cry}

We will next prove a bound of the effective central charge 
$\tilde{c}=c\!-\!24h_{min}$, where $h_{min}$ is the smallest 
conformal weight. 

\begin{cry} Let $V$ be a 
$C_2$-cofinite VOA. 
Then  $\tilde{c}\leq \frac{\dim(V/C_2(V))\!-\!1}{2}$. 
\end{cry}

\pr
The proof is essentially the same as in \cite{GN} with 
slight modifications. Set $k=\dim V/C_2(V)\!-\!1$ and 
define $f_2(q)=\sqrt{2}q^{1/24}\prod_{n=1}^{\infty}(1+q^n)$. 
By using a spanning set of irreducible module $W$ with 
a conformal weight $r$ given in Lemma 2.4, 
there is a polynomial 
$g(q)=\sum_{i=0}^s g_iq^i\in \C[q]$ such that $g_i\geq 0$ and 
$$ {\rm ch}_W(\tau) \leq 2^{-k/2}q^{-k/24}f_2(q)^k g(q)q^r. $$
Here and in the following we shall always 
assume that $0<q<1$. As we showed, 
$${\rm ch}_W(\!-\!1/\tau)=\sum_X a^W_X(\tau){\rm ch}_X(\tau),$$ 
where ${\rm ch}_X(\tau)$ runs over the set of distinct characters 
and $a^W_X(\tau)\in \C[\tau]$. Hence 
$$|\sum a_X(\tau){\rm ch}_X(\tau)|\leq \tilde{q}^{-(k+c)/24}
2^{-k/2}f_2(\tilde{q})^kg(\tilde{q}),$$ 
where $\tilde{q}=e^{-2\pi i/\tau}$ and $f_4(q)
=q^{-1/48}\prod_{n=1}^{\infty}(1\!-\!q^{n-1/2})$. 
In the limit $\tau \to i\infty$ $(q\to 0, \tilde{q}\to 1)$, 
$$ {\rm ch}_W(\tilde{q})=|\tau|^mq^{h-c/24}(a+o(1))g(1) $$
for some integer $m$ and constants $a$, 
where $h$ is a minimal one among conformal weights which appear in 
$\sum_X a^W_X(\tau){\rm ch}_X(\tau)$. Since $\tau \to \!-\!1/\tau$ 
is an involution, there is an irreducible $V$-module $W$ such that 
a character with a minimal conformal weight $h_{min}$ appears in 
${\rm ch}_W(\tilde{q})$. Hence there is a constant $C$ such that 
$$|\tau|^mq^{h_{min}-c/24}|\leq  q^{-k/48}(C+O(q))$$
and so we have $h_{min}\!-\!c/24 \geq \!-\!k/48$ as desired. 
\prend

\end{document}